\documentclass[12pt]{article}

\usepackage{amscd, amssymb, amsthm, amsmath,eucal, verbatim}
\usepackage[dvips]{graphics}

\newcommand{\proj}{\mathbf P}

\newcommand{\rarr}{\rightarrow}
\newcommand{\oh}{{\mathcal{O}}}
\newcommand{\com}{\mathbb{C}}
\newcommand{\Q}{\mathbb{Q}}
\newcommand{\cP}{\mathcal{P}}

\newcommand{\cs}{{\varsigma}}
\newcommand{\Z}{\mathbb{Z}}

\newcommand{\M}{\overline{M}}

\newcommand{\cS}{\mathcal{S}}
\newcommand{\zz}{{\mathfrak{z}}}

\newcommand{\bG}{\mathsf{G}}
\newcommand{\bH}{\mathsf{H}}

\newcommand{\lam}{\lambda}

\newcommand{\A}{\mathcal{A}}
\newcommand{\bA}{\mathsf{A}}
\newcommand{\bAt}{\widetilde{\mathsf{A}}}

\newcommand{\F}{\mathcal{F}}
\newcommand{\cE}{\mathcal{E}}

\newcommand{\Hc}{\bH^\circ}
\newcommand{\Gc}{\bG^\circ}
\newcommand{\al}{\alpha}
\newcommand{\bD}{\mathsf{D}}
\newcommand{\vac}{v_\emptyset}

\newcommand{\lang}{\left\langle}
\newcommand{\rang}{\right\rangle}

\newcommand{\bC}{\mathsf{C}}
\newcommand{\bJ}{\mathsf{J}}
\newcommand{\Pv}{P_\emptyset}
\newcommand{\gli}{\mathfrak{gl}(\infty)}
\newcommand{\nr}[1]{:\!#1\!:}
\newcommand{\ul}{\underline}
\newcommand{\LV}{\Lambda^{\frac\infty2}V}
\newcommand{\LVc}{\Lambda^{\frac\infty2}_0V}
\newcommand{\sh}{{\textstyle \frac12}}
\newcommand{\de}{\delta}
\newcommand{\la}{\lambda}

\newcommand{\MMM}{\mathsf{M}}

\newcommand{\bz}{\mathbf{0}}
\newcommand{\bi}{\boldsymbol{\infty}}

\newcommand{\Ff}[3]{{}_2F_1\left(
    \begin{matrix}
      #1, 1 \\ #2
    \end{matrix} \, ; #3\right)}

\DeclareMathOperator{\ev}{ev} \DeclareMathOperator{\Aut}{Aut}
 \DeclareMathOperator{\Res}{Res}

\newtheorem{Theorem}{Theorem}
\newcommand{\C}{\mathbb{C}}
\newtheorem{Lemma}{Lemma}
\newtheorem{Corollary}[Lemma]{Corollary}
\newtheorem{Proposition}[Lemma]{Proposition}

\DeclareMathOperator{\End}{End}

\numberwithin{equation}{section}

\begin{document}
\title{The equivariant Gromov-Witten theory of $\proj^1$}
\author{A.~Okounkov and R.~Pandharipande}
\date{July 2002}
\maketitle

\setcounter{tocdepth}{2}
\tableofcontents




\setcounter{section}{-1}
\section{Introduction}

\subsection{Overview}

\subsubsection{}

We present here the second in a sequence of three papers devoted
to the Gromov-Witten theory of nonsingular target curves $X$.
Let $\omega \in H^2(X,\Q)$ denote the Poincar\'e dual of the
point class. 
In the
first paper \cite{OP2}, we considered the stationary sector of the
Gromov-Witten theory of $X$ formed by the
descendents of $\omega$. 
The stationary sector was identified in \cite{OP2}
with the Hurwitz theory of $X$ with completed cycles insertions.

The target $\proj^1$ plays a distinguished role in the
Gromov-Witten theory of target curves. Since $\proj^1$ admits a
$\C^*$-action, equivariant localization may be used to study 
Gromov-Witten invariants \cite{GP}. The equivariant Poincar\'e duals,
$$
\bz,\bi\in H^2_{\C^*}(\proj^1, \Q),
$$
of the $\C^*$-fixed points $0,\infty\in\proj^1$
form a basis of the localized
equivariant cohomology of $\proj^1$.  Therefore, the
full equivariant Gromov-Witten theory of $\proj^1$ is quite
similar in spirit to the stationary non-equivariant theory.
Via the non-equivariant limit, the full non-equivariant
theory of $\proj^1$ is captured by the equivariant
theory.

The equivariant Gromov-Witten theory of $\proj^1$ is the subject
of the present paper. We find explicit formulas and establish
connections to integrable hierarchies. The full Gromov-Witten
theory of higher genus target curves will be considered in the
third paper \cite{OP3}. The equivariant theory of $\proj^1$ will
play a crucial role in the derivation of the Virasoro constraints
for target curves in \cite{OP3}.

\subsubsection{}

Our main result here is an explicit operator
description of the equivariant Gromov-Witten theory of $\proj^1$.
We identify all equivariant Gromov-Witten invariants
of $\proj^1$ as vacuum matrix elements
of explicit operators acting in the Fock space
(in the infinite wedge realization).

The result is obtained by combining
the equivariant localization
formula with an  operator formalism for the Hodge integrals
which arise as vertex terms. The
operator formalism for Hodge integrals relies crucially upon
a formula
due to Ekedahl, Lando, Shapiro, and Vainstein (see \cite{ELSV,FP,GV}
and also \cite{OP1})
expressing  basic Hurwitz numbers as Hodge integrals.

\subsubsection{}

As a direct and fundamental consequence of the operator
formalism, we find an integrable hierarchy
governs the equivariant Gromov-Witten theory
of $\proj^1$ ---
specifically,
the $2$--Toda hierarchy of Ueno and Takasaki \cite{UT}.
 The equations of the hierarchy, together with the
string and divisor equations, uniquely determine the
entire theory.

A Toda hierarchy for the non-equivariant
Gromov-Witten of $\proj^1$ was proposed in the
mid 1990's in a series of papers
by the physicists T.~Eguchi, K.~Hori, C.-S.~Xiong,
Y.~Yamada, and S.-K.~Yang on the basis of a conjectural
matrix model description of the theory, see \cite{EgHY,EgY}.
The Toda conjecture was further studied in \cite{P,O,G1,G2} and, for
the stationary sector, proved in \cite{OP2}.

The 2--Toda hierarchy for the equivariant Gromov-Witten theory of
$\proj^1$ obtained
here is both more general and, arguably, more simple than the
hierarchy obtained in the
non-equivariant limit.

\subsubsection{}

The 2--Toda hierarchy governs the equivariant theory of $\proj^1$
just as Witten's KdV hierarchy \cite{W} governs  the Gromov-Witten
theory of a point. However, while the known derivations of the KdV
equations for the point require the analysis of elaborate
auxiliary constructions (see \cite{dF,IZ,K1,O2,OP1}), the Toda
equations for  $\proj^1$ follow directly, almost in textbook
fashion, from the operator description of the theory.

In fact, the Gromov-Witten theory of $\proj^1$ may be viewed as a
more fundamental object than the Gromov-Witten theory of a point.
Indeed, the theory of $\proj^1$ has a simpler and more explicit
structure.
The theory of $\proj^1$  is {\em not} based on the theory of a
point. Rather, the point theory  is perhaps best understood as a
certain special large degree limit case of the $\proj^1$ theory,
see \cite{OP1}.

\subsubsection{}

The proof of the Gromov-Witten/Hurwitz correspondence in \cite{OP2}
assumed a restricted case of the full result: the GW/H correspondence for the
absolute stationary non-equivariant
Gromov-Witten theory of $\proj^1$.
The required case is established here as a direct
consequence of our operator formalism for the equivariant theory of
$\proj^1$ --- completing the proof of the full GW/H correspondence.

While the present paper does not rely upon the
results of \cite{OP2}, much of the motivation can
be found in the study of the stationary theory developed there.

\subsubsection{}

We do not know whether the Gromov-Witten theories of higher genus
target curves are governed by
integrable hierarchies. However, there exist
conjectural Virasoro constraints for the
Gromov-Witten theory of an arbitrary nonsingular projective variety $X$
formulated in 1997
 by Eguchi, Hori, and Xiong (using also ideas of S. Katz), see
\cite{EgHX}.

The results of the present paper will be used in
\cite{OP3} to prove the Virasoro constraints for
 nonsingular target curves $X$.
Givental has recently announced a proof
of  the Virasoro constraints for the
projective spaces $\proj^n$.
These two families of varieties both start with $\proj^1$
but are quite different in flavor.
Curves are of dimension 1, but have non-$(p,p)$ cohomology,
non-semisimple quantum cohomology,
and do not, in general, carry torus actions.
Projective spaces cover all target dimensions, but have
algebraic cohomology, semisimple
quantum cohomology, and always carry torus actions.
Together, these results provide substantial evidence for the
Virasoro constraints.

\subsection{The equivariant Gromov-Witten theory of $\proj^1$}

\subsubsection{}

Let $V= \com \oplus \com$.
Let the algebraic torus $\com^*$ act on $V$ with weights $(0,1)$:
$$
\xi \cdot(v_1,v_2)= (v_1, \xi \cdot v_2)\,.
$$
Let $\proj^1$ denote the projectivization $\proj(V)$.
There is a canonically induced
$\com^*$-action on $\proj^1$.

The $\com^*$-equivariant cohomology ring of a point is $\Q[t]$ where
$t$ is the first Chern class of the standard representation.
The $\com^*$-equivariant cohomology ring $H^*_{\com^*}(\proj^1, \Q)$
is canonically a $\Q[t]$-module.

The line bundle $\oh_{\proj^1}(1)$ admits a canonical $\com^*$-action
which identifies the representation $H^0(\proj^1, \oh_{\proj^1}(1))$ with
$V^*$. Let $h\in H^2_{\com^*}(\proj^1,\Q)$ denote the equivariant first Chern
class of $\oh_{\proj^1}(1)$.
The equivariant cohomology ring of $\proj^1$ is easily determined:
$$H_{\com^*}^*(\proj^1, \Q) = \Q[h, t]/ (h^2 + t h).$$
A free $\Q[t]$-module basis is provided by $1, h$.

\subsubsection{}

Let $\overline{M}_{g,n}(\proj^1,d)$ denote the moduli space of
genus $g$, $n$-pointed stable maps (with connected domains)  to
$\proj^1$ of degree $d$. A canonical $\com^*$-action on
$\overline{M}_{g,n}(\proj^1,d)$ is obtained by translating maps.
The virtual class is canonically defined in   equivariant
homology:
$$[\overline{M}_{g,n}(\proj^1,d)]^{vir} \in H^{\com^*}_{2(2g+2d-2+n)}(
\overline{M}_{g,n}(\proj^1,d), \Q),$$
 where $2g+2d-2+n$ is the expected complex dimension (see, for example,
\cite{GP}).

The equivariant Gromov-Witten theory of $\proj^1$
concerns equivariant integration
over the moduli space   $\overline{M}_{g,n}(\proj^1,d)$.
Two types of equivariant cohomology classes
are integrated.
The {\em primary classes} are:
$$\text{ev}_i^*(\gamma) \in
H^*_{\com^*}(\overline{M}_{g,n}(\proj^1,d), {\mathbb{Q}}),$$
where $\text{ev}_i$ is the morphism defined by evaluation at the
$i$th marked point,
$$
\text{ev}_i: \overline{M}_{g,n}
(\proj^1,d)\rarr \proj^1\,,
$$
and $\gamma\in H^*_{\com^*}(\proj^1, {\mathbb{Q}})$.
The {\em descendent
classes} are:
$$ \psi_i^k \text{ev}_i^*(\gamma),$$
where
$\psi_i\in H^2_{\com^*}(\overline{M}_{g,n}(X,d), {\mathbb{Q}})$
is
the first  Chern class of the cotangent line bundle $L_i$ on
the moduli space of maps.

Equivariant integrals of descendent classes are expressed by
brackets of $\tau_k(\gamma)$ insertions:
\begin{equation}
\label{llqq}
\lang \prod_{i=1}^n \tau_{k_i}(\gamma_i)
\rang_{g,d}^{\circ} =
\int_{[\overline{M}_{g,n}(\proj^1,d)]^{vir}}
\prod_{i=1}^n \psi_{i}^{k_i} \,
\text{ev}_i^*(\gamma_i)\, ,
\end{equation}
where $\gamma_i\in H^*_{\com^*}(\proj^1,\Q)$.
As in \cite{OP2}, the superscript $\circ$
indicates the connected theory. The theory
with possibly disconnected domains is denoted
by $\lang\,\, \rang^\bullet$. The
equivariant integral in (\ref{llqq}) denotes equivariant push-forward to
a point. Hence,
the bracket takes values in
$\Q[t]$.

\subsubsection{}

We now define the equivariant Gromov-Witten potential $F$ of
$\proj^1$. Let $z,y$ denote the variable sets,
$$
\{z_0,z_1,z_2, \ldots\}, \ \ \{y_0, y_1, y_2, \ldots \}.
$$
The variables $z_k$, $y_k$ correspond to the descendent insertions
$\tau_k(1)$, $\tau_k(h)$ respectively.
Let $T$ denote the formal sum,
$$
T= \sum_{k=0}^\infty z_k \tau_k(1)+ y_k \tau_k(h)\,.
$$
The potential is a generating series of equivariant integrals:
$$
F= \sum_{g=0}^\infty  \sum_{d=0}^\infty \sum_{n=0}^\infty
u^{2g-2} q^d \lang \frac{T^n}{n!} \rang_{g,d}^{\circ}\,.
$$
The potential $F$ is an element of $\Q[t][[z,y,u,q]]$.

\subsubsection{}

The (localized) equivariant
cohomology of $\proj^1$ has a canonical basis provided by the classes,
$$
\bz,\bi\in H^2_{\C^*}(\proj^1)\,,
$$
of
Poincar\'e duals of the $\C^*$-fixed points $0,\infty\in\proj^1$.
An elementary calculation yields:
\begin{equation} \label{nv}
\bz = t \cdot 1 + h, \ \ \bi = h.
\end{equation}

Let $x_i$, $x^\star_i$ be the variables corresponding to
the descendent insertions $\tau_k(\bz)$, $\tau_k(\bi)$, respectively. The
variable sets $x, x^\star$ and $z,y$ are related by the transform dual to \eqref{nv},
$$x_i = \frac{1}{t} z_i , \ \ x_i^\star= -\frac{1}{t}z_i +y_i.$$
The equivariant Gromov-Witten potential  of
$\proj^1$ may be written in the $x_i$, $x_i^\star$ variables as:
$$
F= \sum_{g=0}^\infty  \sum_{d=0}^\infty
u^{2g-2} q^d \lang
\exp\left(
\sum_{k=0}^\infty x_k \tau_k(\bz)+ x^\star_k \tau_k(\bi)
\right)\rang_{g,d}^{\circ}\,.
$$

\subsection{The equivariant Toda equation}

\subsubsection{}

Let the classical series $F^c$ be the genus 0, degree 0,
3-point summand of
$F$ (omitting $u,q$).
The
classical series generates the
equivariant integrals of triple products
in $H^*_{\com^*}(\proj^1, \Q)$. We find,
$$
F^{c}=
\frac{1}{2}z_0^2 y_0 - \frac{1}{2}t z_0 y_0^2 + \frac{1}{6}t^2 y_0^3\,.
$$
The classical series  does not depend upon $z_{k>0}$, $y_{k>0}$.

Let $F^0$ be the genus 0 summand of $F$ (omitting
$u$).
The {\em small phase space} is the hypersurface defined by the
conditions:
$$
z_{k>0}=0,\ y_{k>0}=0\,.
$$
The restriction of the genus 0 series to the
small phase space is easily calculated:
$$
F^0\big|_{z_{k>0}=0,\ y_{k>0}=0} = F^{c}+  qe^{y_0}\,.
$$
The second derivatives of the restricted function $F^0$ are:
$$
F^0_{z_0z_0} = y_0, \ \ F^0_{z_0y_0} = z_0 - t y_0, \ \
F^0_{y_0y_0} = -t z_0 + t^2 y_0 +qe^y_0.$$ Hence, we find  the
equation
\begin{equation}
\label{gztd}
t F^0_{z_0y_0} + F^0_{y_0y_0}  = q\exp( F^0_{z_0z_0} )
\end{equation}
is valid at least on the small phase space.

\subsubsection{}

The equivariant Toda equation for the full equivariant potential $F$
takes
a similar form:
\begin{equation}
\label{tda}
t F_{z_0y_0} + F_{y_0y_0} = \frac{q}{u^{2}}  \exp( F(z_0+u)+ F(z_0-u) - 2F ),
\end{equation}
where $F(z_0\pm u)= F(z_0\pm u, z_1, z_2, \ldots, y_0, y_1, y_2, \ldots,u,q).$
In fact, the equivariant Toda equation
specializes to (\ref{gztd}) when restricted to genus 0 and the small
phase space.

\subsubsection{}

In the variables $x_i$, $x^\star_i$, the equivariant Toda equation
may be written as:
\begin{equation}
  \label{eqT}
  \frac{\partial^2}{\partial x_0\, \partial x^\star_0} \, F =
\frac{q}{u^2} \, \exp\left(\Delta F\right)\,.
\end{equation}
Here, $\Delta$ is the difference operator,
$$
\Delta  = e^{u\partial} - 2 + e^{-u\partial} \,,
$$
and
$$
\partial = \frac{\partial}{\partial z_0} =
\frac1{t} \, \left(\frac{\partial}{\partial x_0} -
\frac{\partial}{\partial x^\star_0}\right)
$$
is the vector field creating a $\tau_0(1)$ insertion.

The equivariant Toda equation in form \eqref{eqT}
is recognized as the 2--Toda equation: obtained from the
standard Toda equation by replacing the second time derivative
by $\frac{\partial^2}{\partial x_0\, \partial x^\star_0}$.
The $2$--Toda equation is a 2-dimensional time analogue 
 of the standard Toda equation.

\subsubsection{}

A central result of the paper is the derivation of the $2$--Toda equation for
the equivariant theory of $\proj^1$.

\vspace{+10pt}
\noindent{\bf Theorem.}
 The equivariant Gromov-Witten potential  of $\proj^1$ satisfies the
2--Toda equation \eqref{eqT}.
\vspace{+10pt}

The 2--Toda equation is a strong constraint. Together with the equivariant divisor
and string equations, the 2--Toda determines $F$ from the degree 0 invariants,
see \cite{P}.

The 2--Toda equation arises as the lowest equation in a
hierarchy of partial differential equations
identified with the 2--Toda hierarchy of Ueno and Takasaki \cite{UT},
see Theorem \ref{TT2} in Section \ref{fififi}.

\subsection{Operator formalism}

\subsubsection{}

The 2--Toda equation \eqref{eqT} is a direct consequence of
the following operator formula for the equivariant Gromov-Witten
theory of $\proj^1$:
\begin{equation}
  \label{eF}
\exp F
 =  \lang e^{\sum x_i \bA_i} \,
\, e^{\al_1} \, \left(\frac{q}{u^2}\right)^H \, e^{\al_{-1}} \,
e^{\sum x^\star_i \bA^\star_i}
\rang \,.
\end{equation}
Here,
$\bA_i$, $\bA^\star_i$, and $H$ are explicit
operators in the Fock space.
The brackets $\lang \, \rang $ denote
the vacuum matrix element. The operators $\bA$, which
depend on the parameters $u$ and $t$, are constructed
in Sections \ref{ttt} and  \ref{ththth}. The exponential
$e^F$ of the equivariant potential is called the $\tau$-function
of the theory.

The operator formula \eqref{eF}, stated as Theorem
\ref{vvv} in Section \ref{ththth},
is fundamentally the main result of the paper.

\subsubsection{}

In our previous paper \cite{OP2}, the stationary non-equivariant
Gromov-Witten theory of $\proj^1$ was expressed as a
similar vacuum expectation.
The equivariant formula \eqref{eF} specializes to the
absolute case of the operator formula of \cite{OP2} when the equivariant
parameter $t$ is set to zero. Hence, the equivariant formula \eqref{eF}
completes the proof of the Gromov-Witten/Hurwitz
correspondence discussed in \cite{OP2}.

\subsection{Plan of the paper}

\subsubsection{}

In Section \ref{ooo}, the virtual localization formula
of \cite{GP} is applied to express the
equivariant $n+m$-point function
as a graph sum with vertex Hodge integrals. Since $\proj^1$
has two fixed points, the graph sum reduces to a sum over partitions.

Next, an operator formula for Hodge integrals is obtained in
Section \ref{ttt}. A starting point here is provided by
the Ekedahl-Lando-Shapiro-Vainstein formula expressing
the necessary Hodge integrals as Hurwitz numbers.
The main result of the section is Theorem \ref{hodfor}
which expresses the generating function for Hodge integrals
as a vacuum matrix element of a product of
explicit operators $\A$ acting on the infinite wedge
space.

Commutation relations for
the  operators $\A$ are required in the proof of Theorem \ref{hodfor}.
The technical derivation of these commutation relations is postponed
to  Section \ref{ffff}.

In Section \ref{ththth},
the operator formula for Hodge integrals is  combined
with the results of Section \ref{ooo}
to obtain Theorem \ref{vvv}, the operator formula
for the equivariant Gromov-Witten theory of $\proj^1$.

The 2--Toda equation \eqref{eqT} and the full 2--Toda
hierarchy 
are deduced from Theorem \ref{vvv}
in Section \ref{fififi}.


\subsubsection{}

We follow the notational conventions of
\cite{OP2} with one important difference. The letter
$\bH$ is used here to denote the generating function
for Hodge integral, whereas $\bH$  was
used to denote Hurwitz numbers in \cite{OP2}.

\subsection{Acknowledgments}

We thank E.~Getzler and A.~Givental for discussions of the
Gromov-Witten theory of $\proj^1$.
In particular, the explicit form
of the linear change of time variables appearing
in the equations of the 2--Toda hierarchy (see 
Theorem \ref{TT2}) was previously conjectured by
Getzler in \cite{G2}.

A.O.\ was partially supported by
DMS-0096246 and fellowships from the Sloan and Packard foundations.
R.P.\ was partially supported by DMS-0071473
and fellowships from the Sloan and Packard foundations.

The paper was completed during a visit to the Max Planck
Institute in Bonn in the summer of 2002.



\section{Localization for $\proj^1$}
\label{ooo}
\subsection{Hodge integrals}

\subsubsection{}

Hodge integrals of the $\psi$ and $\lambda$ classes
over the moduli space of curves arise as
vertex terms in the localization formula for Gromov-Witten
invariants of $\proj^1$.

Let $L_i$ be the $i${th} cotangent line bundle on $\M_{g,n}$.
The $\psi$ classes are defined by:
$$
\psi_i = c_1(L_i) \in H^2(\M_{g,n}, \Q)\,.
$$
Let $\pi: C \rarr \overline{M}_{g,n}$ be the universal curve.
Let $\omega_\pi$ be the relative dualizing sheaf.
Let ${\mathbb E}$ be the rank $g$ Hodge bundle on the moduli space
$\overline{M}_{g,n}$,
$${\mathbb E} = \pi_*(\omega_\pi).$$
The $\lambda$ classes are defined by:
$$\lambda_i = c_i({\mathbb E})\in H^*(\overline{M}_{g,n}, \Q).$$

Only Hodge integrands linear in the $\lambda$ classes arise
in the localization formula for $\proj^1$.
Let $\Hc_g(z_1, \ldots, z_n)$ be the $n$-point function of $\lambda$-linear Hodge integrals over
the moduli space $\M_{g,n}$:
$$
\Hc_g(z_1, \ldots, z_n)
={\textstyle{\prod}} z_i \int_{\M_{g,n}} \frac{1-\la_1+\la_2-\dots\pm \la_g}
{\prod (1- z_i \psi_i)} \,.
$$
Note the shift of indices caused by the product $\prod z_i$.

\subsubsection{}

The function $\Hc_g(z)$ is defined for all $g, n\geq 0$.
Values
corresponding to unstable moduli spaces are set by definition.
All $0$-point functions $\Hc_g()$, both stable and unstable,
vanish. The unstable 1 and 2-point functions are:
\begin{equation}
\Hc_0(z_1)=\frac{1}{z_1}\,, \quad
\Hc_0(z_1,z_2)=\frac{z_1 z_2}{z_1+z_2}\,.\label{hand}
\end{equation}

\subsubsection{}

Let $\Hc(z_1, \ldots, z_n,u)$ be the full $n$-point function 
of $\lambda$-linear Hodge integrals:
$$
\Hc(z_1,\dots,z_n,u)=\sum_{g\ge 0} u^{2g-2} \, \Hc_g(z_1, \ldots, z_n)\,.
$$
Let $\bH(z_1, \ldots, z_n,u)$ be the corresponding disconnected 
$n$-point function.
The disconnected $0$-point function is defined by:
$$\bH(u)= 1,$$
For $n>0$, the disconnected $n$-point function is defined by:
$$
\bH(z_1, \ldots, z_n,u)= 
\sum_{P\in \text{Part}[n]} \prod_{i=1}^{\ell(P)}
\Hc(z_{P_i},u),$$
where $\text{Part}[n]$ is the set of partitions $P$
 of the set $\{1, \ldots,n\}$. Here, $\ell(P)$ is the length of the
partition, and $z_{P_i}$ denotes the variable
set indexed by the part $P_i$.
The genus expansion for the disconnected function,
\begin{equation}
\bH(z_1,\dots,z_n,u)=
\sum_{g\in \Z} u^{2g-2} \, \bH_g(z_1, \ldots, z_n)\,,
\label{bHu}
\end{equation}
contains negative genus terms.

\subsection{Equivariant $n+m$-point functions}

\subsubsection{}

Let $\Gc_{g,d}(z_1, \ldots, z_n, w_1, \ldots, w_m)$ be the
 $n+m$-point function of genus $g$, degree $d$ equivariant
Gromov-Witten invariants of $\proj^1$ in the basis determined by
$\bz$ and $\bi$:
$$
\Gc_{g,d}(z,w)= {\textstyle{\prod}} z_i
{\textstyle{\prod}} w_j
 \int_{[\M_{g,n+m}(\proj^1,d)]^{vir}}
\prod \frac{\ev_i^*(\bz)}{1-z_i\,\psi_i} \,
\prod \frac{\ev_j^*(\bi)}{1-w_j\,\psi_j}  \,.
$$
The values corresponding to unstable moduli spaces are set by definition.
The unstable $0$-point functions are set to 0:
\begin{equation}
\Gc_{0,0}()=0\, , \quad
\Gc_{1,0}()=0 \, . \label{gunst0}
\end{equation}
The unstable $1$ and $2$-point functions are:
\begin{alignat}{2}
\label{gunst}
\Gc_{0,0}(z_1)&=\frac{1}{z_1}\,,& \quad
\Gc_{0,0}(w_1)&=\frac{1}{w_1}\, , 
\end{alignat}
\begin{alignat}{3}
\Gc_{0,0}(z_1,z_2)&=\frac{tz_1 z_2}{z_1+z_2}\,, &\quad
\Gc_{0,0}(z_1,w_1)&=0\,, &\quad
\Gc_{0,0}(w_1,w_2)&=\frac{tw_1 w_2}{w_1+w_2}\,. \notag
\end{alignat}
These values
will be seen to be compatible with the special
values \eqref{hand}.

\subsubsection{}

The $n+m$-point function $\Gc_{g,d}(z,w)$ is defined for all $g,d,n,m\geq 0$.
The $0$-point function $\Gc_{0,1}()$ is nontrivial since
$$\Gc_{0,1}() =\langle \rangle_{0,1}^\circ =1 .$$
In fact, $\Gc_{0,1}()$ is the {\em only} nonvanishing 0-point function
for $\proj^1$.

Let $\Gc_d(z,w,u)$ be the full $n+m$-point function for equivariant degree $d$
Gromov-Witten invariants $\proj^1$:
$$
\Gc_d(z_1,\dots,z_n,w_1,\dots,w_m,u)=
\sum_{g\ge 0} u^{2g-2} \, \Gc_{g,d} (z_1, \ldots, z_n,w_1, \ldots, w_m)\, .
$$
The only nonvanishing $0$-point functions is:
$$\Gc_1() = u^{-2}.$$

\subsubsection{}

Let $\bG_d(z,w,u)$ be the corresponding disconnected
$n+m$-point function. The degree 0, $0$-pointed disconnected
function is defined by:
$$\bG_0(u)=1.$$
In all other cases,
$$
\bG_d(z_1, \ldots, z_n, w_1, \ldots, w_m,u)=
\sum_{P\in \text{Part}_d[n,m]} \frac{1}{|\Aut(P)|} \prod_{i=1}^{\ell(P)}
\Gc_{d_i}(z_{P_i},w_{P'_i},u).$$
An element $P\in \text{Part}_d[n,m]$ consists of the data
$$\{ (d_1,P_1,P'_1)  \ldots, (d_\ell, P_\ell,P'_\ell)\}\, ,$$
where $d_i$ is a non-negative degree partition,
$$\sum_{i=1}^l d_i =d,$$
and $\{P_i\}$ and $\{P'_i\}$ are set partitions
with the empty set as an allowed part,
$$
\bigcup_{i=1}^l P_i= \{1, \ldots, n\}, \ \
\bigcup_{i=1}^l P'_i = \{1, \ldots, m\}\,.
$$

\subsubsection{}

Two remarks about the $n+m$-point function $\bG_d(z,w,u)$ are in order.
First, $\bG_d$ systematically includes the
unstable contributions (\ref{gunst}). These contributions
will later have to be
removed to study the true equivariant Gromov-Witten theory.
However, the inclusion of the unstable contributions here
will simplify many formulas.
Second,
the $0$-point function $\Gc_{1}()$ contributes to all disconnected
functions $\bG_d$ for positive $d$. For example:
$$
\bG_{2}(z_1)=\Gc_2(z_1)+ \Gc_1(z_1) \, \Gc_1() + \Gc_0(z_1)\,
\frac{\Gc_1()^2}{2}\,.
$$
These occurrences of $\Gc_1()$ provide no difficulty.

\subsection{Localization: vertex contributions}

\subsubsection{}

The localization formula for $\proj^1$
expresses the $n+m$-point function $\bG_d(z,w,u)$
as an automorphism-weighted sum over bipartite graphs
with vertex Hodge integrals. We refer the reader to \cite{GP} for a discussion
of localization in the context of virtual classes.
The localization formula for $\proj^1$ is explicitly treated
in \cite{GP,OP1}.

\subsubsection{}

Let $\Gamma$ be a graph arising in the localization formula
for the virtual class $[\M_{g,n+m}(\proj^1,d)]^{vir}$.
Let $v_0$ be a vertex of $\Gamma$
 lying over the fixed point $0\in \proj^1$.
We will study  the vertex contribution $C(v_0)$ to the
equivariant integral
\begin{equation}
\label{pqwr}
\prod z_i \prod w_j
 \int_{[\M_{g,n+m}(\proj^1,d)]^{vir}}
\prod \frac{\ev_i^*(\bz)}{1-z_i\,\psi_i} \,
\prod \frac{\ev_j^*(\bi)}{1-w_j\,\psi_j}  \,.
\end{equation}
For a vertex $v_\infty$ lying over $\infty \in \proj^1$,
the vertex contribution $C(v_\infty)$
is obtained simply by exchanging the
roles of $z$ and $w$ and applying the transformation
$t\mapsto -t$.

Each vertex $v_0$ of the localization graph $\Gamma$ carries several  additional structures:
\begin{enumerate}
\item[$\bullet$] $g(v_0)$, a genus assignment,
\item[$\bullet$] $e(v_0)$ incident edges of degrees $d_1, \ldots, d_{e(v_0)}$,
\item[$\bullet$] $n(v_0)$ marked points indexed by $I(v_0) \subset \{1, \ldots, n\}$.
\end{enumerate}
The data contribute factors to the vertex contribution $C(v_0)$ according to
the following table:

{
\renewcommand{\arraystretch}{3}
\begin{center}
\begin{tabular}{|l|l|}
\hline
$\displaystyle
t^{g(v_0)-1} \left( \sum_{i=1}^{g(v_0)} (-1)^i \frac{\la_i}{t^i}
\right)$ &
determined by the genus $g(v_0)$ \\
$\displaystyle
\frac{d_i^{d_i}\, t^{-d_i}}{d_i!}\,
\frac{t d_i}{t- d_i \, \psi_i}$
& for each edge of degree $d_i$ \\
$\displaystyle
\frac{t z_i}{1- z_i \, \psi_i}$
& for each marking $i\in I(v_0)$ \\
\hline
\end{tabular}
\end{center}
}

The vertex contribution $C(v_0)$ is obtained by
multiplying the above factors and integrating over the moduli space
$\M_{g(v_0),val(v_0)}$ where $$val(v_0)=e(v_0)+n(v_0).$$

\subsubsection{}

By the dimension constraint for the integrand,
$$
\dim \ \M_{g(v_0),val(v_0)} = 3g(v_0)-3+val(v_0)\,,
$$
the vertex integral is unchanged by
the transformation
$$
\psi_i \mapsto t \psi_i \,, \quad \la_i \mapsto t^i \la_i \,,
$$
together with a division by $t^{3g(v_0)-3+val(v_0)}$.
The vertex contribution $C(v_0)$ then takes the following form:
\begin{multline*}
  \frac{
\prod_{i=1}^{e(v_0)} d_i^{d_i}\big/d_i !}
{\displaystyle
t^{2 g(v_0)-2+ d(v_0)+val(v_0)}} \times\\
\int_{\M_{g(v_0),val(v_0)}}
  \left( \sum_{i=1}^{g(v_0)} (-1)^i {\la_i}{} \right)
  \prod_{i=1}^{e(v_0)}\frac{d_i}{1-d_i \psi_i} \prod_{i\in I(v_0)}
  \frac{t z_i}{1- t z_i \psi_i},
\end{multline*}
where $d(v_0)=\sum_{i=1}^{e(v_0)} d_i$ is the total degree of $v_0$.
We may rewrite $C(v_0)$ in terms of $\Hc_{g(v_0)}$:
\begin{equation}
\label{paap}
C(v_0)=  \frac{
\prod_{i=1}^{e(v_0)} d_i^{d_i}\big/d_i !}
{\displaystyle
t^{2 g(v_0)-2+ d(v_0)+val(v_0)}} \,
\Hc_{g(v_0)}(d_1,\dots,d_{e(v_0)},\dots,tz_i,\dots).
\end{equation}
Since the $val(v_0)$-point function $\Hc_{g(v_0)}$ is defined for all
$g(v_0),val(v_0) \geq 0$, we can define the vertex contribution $C(v_0)$
by (\ref{paap})
in case the moduli space
$\M_{g(v_0), val(v_0)}$ is unstable. We note $C(v_0)$ vanishes if $val(v_0)=0$.

\subsection{Localization: global formulas}

\subsubsection{}

Let $\Gamma$ be a graph arising in the localization formula for $[\M_{g,n+m}(\proj^1,d)]^{vir}$. Let
$$V(\Gamma)= V_0(\Gamma) \cup V_\infty(\Gamma)$$ be the vertex set
divided by fixed point assignment. Let $E(\Gamma)$ be the edge set. Let $d_e$ be the
degree of an edge $e$.
The graph $\Gamma$ satisfies
 three global properties:
\begin{enumerate}
\item[$\bullet$] a genus condition,
$\sum_{v\in V(\Gamma)} (2g(v)-2+e(v)) =2g-2$,
\item[$\bullet$] a degree condition,
 $\sum_{v\in V(\Gamma)} d(v) = 2d$,
\item[$\bullet$] a marking condition, $\bigcup_{v_0\in V_0(\Gamma)} I(v_0)= \{1, \ldots,n\}$
\ \ (similarly for $\infty$).
\end{enumerate}
The contribution of $\Gamma$ to the integral (\ref{pqwr}) is:
$$\frac{1}{\prod_{e\in E(\Gamma)} d_e} \
\frac{1}{|{\text {Aut}}(\Gamma)|} \prod_{v\in V(\Gamma)} C(v).$$
As the integral (\ref{pqwr}) is over the moduli space of maps with connected
domains, $\Gamma$ must also be connected. If disconnected domains
are allowed for stable maps, the graphs $\Gamma$ are also allowed to
be disconnected.

\subsubsection{}

The $n+m$-point functions $\bG_d$ may be now expressed in terms of the functions
$\bH$.

\begin{Proposition}
\label{Loc}
For $d\geq 0$,  we have
\begin{multline}\label{bGd}
   \bG_d(z_1,\dots,z_n,w_1,\dots,w_m,u) = \frac{1}{\zz(\mu)}
\times \\
\sum_{|\mu|=d} \frac{
    (u/t)^{\ell(\mu)}\,
(-u/t)^{\ell(\mu)} } {t^{d+n} (-t)^{d+m}
} \, \left(\prod
    \frac{\mu_i^{\mu_i}}{\mu_i !}\right)^2 \bH(\mu,tz,\tfrac ut) \,
  \bH(\mu,-tw,-\tfrac ut)\,.
\end{multline}
\end{Proposition}
\noindent
The summation in \eqref{bGd}
is over all partitions $\mu$ of $d$,
$\ell(\mu)$ denotes the number of parts of $\mu$ and
$$
\zz(\mu)= \left|\Aut(\mu)\right| \prod_{i=1}^{\ell(\mu)} \mu_i
$$
where $\Aut(\mu) \cong \prod_{i\geq 1} S({m_i(\mu)})$
is the symmetry group permuting equal parts of the
partition $\mu$.

\begin{proof} Each degree $d$, possibly disconnected, localization graph $\Gamma$
yields a partition $\mu$ of $d$ obtained from the edge degrees.
The sum over localization graphs with a fixed edge degree partition $\mu$
can be evaluated by the
vertex contribution formula (\ref{paap}) together with
the global graph constraints.
The result is exactly the $\mu$ summand
in \eqref{bGd}
(the edge and graph automorphisms are
incorporated in the prefactors).
The Proposition is then a restatement of the virtual localization formula:
equivariant integration against the virtual class is obtained
by summing over all localization graph contributions.
\end{proof}

The degree $0$ localization formula is special as the graphs are edgeless.
However, with our  conventions regarding $0$-pointed functions,
Proposition \ref{Loc} holds without modification.
We find, for example,
$$
\bG_0(z_1,\dots,z_n,u) = t^{-n} \, \bH(tz,\tfrac ut) \, .
$$
In particular, the definitions of the unstable contributions for $\bG$ and
$\bH$
are compatible.



\section{Operator formula for Hodge integrals}
\label{ttt}



We will express Hodge integrals as matrix elements in the infinite wedge space.
The basic properties of the infinite wedge space and our notational conventions
are summarized in Section \ref{iwr}. A discussion can also be found in Section 2 of
\cite{OP2}.

\setcounter{subsection}{-1}

\subsection{Review of the infinite wedge space}

\label{iwr}

\subsubsection{}

Let $V$ be a linear space with basis $\left\{\ul{k}\right\}$
indexed by the half-integers:
$$V = \bigoplus_{k \in \Z+ \sh} \com \, \ul{k}.$$
For each subset $S=\{s_1>s_2>s_3>\dots\}\subset \Z+\sh$ satisfying:
\begin{enumerate}
\item[(i)] $S_+ = S \setminus \left(\Z_{\le 0} - \sh\right)$ is
finite,
\item[(ii)]$S_- =
\left(\Z_{\le 0} - \sh\right) \setminus S$ is finite,
\end{enumerate}
we denote by  $v_S$ the following infinite wedge product:
\begin{equation}
v_S=\ul{s_1} \wedge \ul{s_2} \wedge  \ul{s_3} \wedge \dots\, .
\label{vS}
\end{equation}
By definition,
 $$\LV= \bigoplus \com \,v_S
$$
is the linear space with basis $\{v_S\}$.
Let $(\,\cdot\,,\,\cdot\, )$ be the inner product on  $\LV$ for which
$\{v_S\}$ is an orthonormal basis.

\subsubsection{}

The fermionic operator $\psi_k$ on $\LV$ is defined by wedge product with
the vector $\ul{k}$,
$$
\psi_k \cdot v = \ul{k} \wedge v  \,.
$$
The operator $\psi_k^*$ is defined as the adjoint of $\psi_k$
with respect to the inner product  $(\,\cdot\,,\,\cdot\, )$.

These operators satisfy the canonical anti-commutation relations:
\begin{gather}
  \psi_i \psi^*_j + \psi^*_i \psi_j = \delta_{ij}\,, \\
  \psi_i \psi_j + \psi_j \psi_1 = \psi_i^* \psi_j^* + \psi_j^*
  \psi_i^*=0.
\end{gather}
The {\em normally
ordered} products are defined by:
\begin{equation}\label{e112}
\nr{\psi_i\, \psi^*_j} =
\begin{cases}
\psi_i\, \psi^*_j\,, & j>0 \,,\\
-\psi^*_j\, \psi_i\,, & j<0 \,.
\end{cases}
\end{equation}

\subsubsection{} \label{sCH}

Let $E_{ij}$, for $i,j\in \Z+\sh$, be the standard basis of
matrix units of $\gli$.
The assignment
$$
E_{ij} \mapsto\, \nr{\psi_i\, \psi^*_j}\ \ ,
$$
defines
a projective representation of the Lie algebra $\gli=\mathfrak{gl}(V)$ on
$\LV$.

The {\em charge} operator $C$ corresponding to the identity matrix of $\gli$,
$$
C= \sum_{k\in\Z+\frac12} \, E_{kk},
$$
acts on the basis $v_S$ by:
$$
C \, v_S = (|S_+| - |S_-|) v_S \,.
$$
The kernel of $C$, the zero charge subspace,
is spanned by the vectors
\begin{equation*}
v_{\lam}=\ul{\lam_1-\tfrac12} \wedge \ul{\lam_2-\tfrac32} \wedge
\ul{\lam_3-\tfrac52} \wedge \dots
\end{equation*}
indexed by all partitions $\la$. We will denote the kernel by $\LVc$.

The eigenvalues on $\LVc$ of the {\em energy} operator,
$$
H= \sum_{k\in\Z+\frac12} k \, E_{kk},
$$
are easily identified:
$$
H \, v_\lambda = |\lambda| \, v_\lambda\,.
$$
The vacuum vector
$$
\vac = \ul{-\tfrac12} \wedge \ul{-\tfrac32} \wedge \ul{-\tfrac52}
\wedge \dots
$$
is the unique vector with the minimal (zero) eigenvalue of $H$.

The {\em vacuum expectation} $\lang A \rang$ of an operator $A$ on $\LV$ is defined
by  the inner product:
$$\lang A \rang = (A v_\emptyset, v_\emptyset).$$

\subsubsection{}

For any $r\in\Z$, we define
\begin{equation}
\cE_r(z) = \sum_{k\in\Z+\frac12} \, e^{z(k-\frac{r}2)} \,  E_{k-r,k}
+ \frac{\delta_{r,0}}{\cs(z)} \, ,\label{Er}
\end{equation}
where the function $\cs(z)$ is defined by
\begin{equation}
  \label{cs}
  \cs(z) = e^{z/2} - e^{-z/2} \,.
\end{equation}
The exponent in \eqref{Er} is set to satisfy:
$$
\cE_r(z)^* = \cE_{-r}(z)\,,
$$
where the adjoint is with respect to the standard inner
product on $\LV$.

Define the operators $\cP_k$ for $k>0$
by:
\begin{equation}
\frac{\cP_k}{k!} = \, [z^k] \, \cE_0(z)\,,\label{defcP}
\end{equation}
where $[z^k]$ stands for the coefficient of $z^k$.
The operator,
$${\mathcal F}_2 = \frac{\cP_2}{2!}= \sum_{k\in\Z+\frac12} \frac{k^2}{2} E_{k,k}\, ,$$
will play a special role.

\subsubsection{}

The operators $\cE$ satisfy the following fundamental
commutation relation:

\begin{equation}\label{commcEE}
  \left[\cE_a(z),\cE_b(w)\right] =
\cs\left(\det  \left[
\begin{smallmatrix}
  a & z \\
b & w
\end{smallmatrix}\right]\right)
\,
\cE_{a+b}(z+w)\,.
\end{equation}

Equation \eqref{commcEE}
automatically incorporates the central
extension of the $\gli$-action, which appears as the
constant term in $\cE_0$ when $r=-s$.

\subsubsection{}

The operators $\cE$ specialize to the
standard bosonic operators on $\LV$:
$$
\al_k = \cE_k(0)\,, \quad k\ne 0 \,.
$$
The commutation relation \eqref{commcE} specializes to
the following equation
\begin{equation}
\label{qhhqq}
[\al_k, \cE_l(z)] =
 \cs(kz) \, \cE_{k+l}(z) \,.
\end{equation}
When $k+l=0$, equation (\ref{qhhqq}) has the following constant term:
$$
\frac{\cs(kz)}{\cs(z)}=\frac{e^{kz/2}-e^{-kz/2}}{e^{z/2}-e^{-z/2}} \,.
$$
Letting $z\to 0$, we recover the standard relation:
\begin{equation*}
[\al_k,\al_l]= k \, \delta_{k+l} \,.
\end{equation*}
%


\subsection{Hurwitz numbers and Hodge integrals}

\subsubsection{}

Let $\mu$ be a partition of size $|\mu|$ and length $\ell(\mu)$. Let
$\mu_1, \ldots, \mu_{\ell}$ be the parts of $\mu$.
Let $\bC_g(\mu)$ be the Hurwitz number of genus $g$, degree $|\mu|$,
covers of $\proj^1$ with
profile $\mu$ over $\infty\in\proj^1$ and
simple ramifications over
$$
b=2g+|\mu|+\ell(\mu)-2
$$
fixed points of ${\mathbf{A}^1}\subset \proj^1$.
The Hurwitz number $\bC_g(\mu)$ is a automorphism-weighted count
of possibly disconnected covers (the genus $g$ may be negative).
The Ekedahl-Lando-Shapiro-Vainstein
formula expresses $\bC_g(\mu)$ in terms of $\lambda$-linear
Hodge integrals:
\begin{equation}
\label{elsv}
\bC_g(\mu) = \frac{b!}{\zz(\mu)}
\, \left(\prod \frac{\mu_i^{\mu_i}}{\mu_i !} \right)\,
\bH_g(\mu_1,\dots, \mu_\ell) \,,
\end{equation}
see \cite{ELSV} or \cite{FP,GV} for a Gromov-Witten theoretic approach.

\subsubsection{}

The Hurwitz numbers $\bC_g(\mu)$ admit a standard expression
in terms of the characters of the symmetric group.
The  character formula may be rewritten as a vacuum expectation in
the infinite wedge space:
\begin{equation}
\label{jkk}
\bC_g(\mu) =
\frac{1}{\zz(\mu)}
\lang e^{\al_1} \F_2^b \prod \al_{-\mu_i} \rang \, .
\end{equation}
A derivation of (\ref{jkk}) can be found, for example, in \cite{O,OP2}.
Using the ELSV formula (\ref{elsv}), we find,
$$
\bH(\mu_1,\dots,\mu_\ell,u) = u^{-|\mu|-\ell(\mu)} \,
\left(\prod \frac{\mu_i !}{\mu_i^{\mu_i}} \right)\,
\lang e^{\al_1} e^{u \F_2} \prod \al_{-\mu_i} \rang \,.
$$

\subsubsection{}

Since the operators $e^{-\al_1}$  and $e^{-u \F_2}$ fix the vacuum
vector, we may rewrite the last equation as:
\begin{multline}
\label{lllst}
\bH(\mu_1,\dots,\mu_\ell,u) = \\
u^{-|\mu|-\ell(\mu)} \,
\left(\prod \frac{\mu_i !}{\mu_i^{\mu_i}} \right)\,
\lang\prod \Big( e^{\al_1} e^{u \F_2}  \al_{-\mu_i} e^{-u \F_2} e^{-\al_1}\Big) \rang \,.
\end{multline}
Equation \eqref{lllst} holds, by construction, for positive integer
values of $\mu_i$. We will rewrite the right side and reinterpret
\eqref{lllst} as an equality of analytic functions of $\mu$.

\subsection{The operators $\A$}

\subsubsection{}

The following operators will play a central role in the paper:
\begin{equation}\label{defA}
\A(a,b) = \cS(b)^a \,
\sum_{k\in \Z} \frac{\cs(b)^{k}}{(a+1)_k} \, \cE_k(b)
 \,,
\end{equation}
where $a$ and $b$ are parameters and
$$
\cs(z)=e^{z/2}-e^{-z/2}\,, \quad \cS(z)=\frac{\cs(z)}{z} =
\frac{\sinh{z/2}}{z/2} \,.
$$
In \eqref{defA}, we use the standard notation:
$$
(a+1)_k = \frac{(a+k)!}{a!} =
\begin{cases}
 (a+1) (a+2) \cdots (a+k) \,, & k\ge 0 \,,  \\
 (a (a-1) \cdots (a+k+1))^{-1}  \,, & k \le 0 \,.
\end{cases}
$$
If $a\ne 0,1,2,\dots$, the sum in \eqref{defA} is
infinite in both directions. If $a$ is a nonnegative
integer, the summands with $k \le - a - 1$ in \eqref{defA}
vanish.

\subsubsection{}
Definition \eqref{defA} is motivated by
the following result.

\begin{Lemma}\label{conjF}
For $m=1,2,3,\dots$, we have
$$
e^{\al_1} \, e^{u \F_2} \, \al_{-m} \, e^{-u \F_2} \, e^{-\al_1} =
\frac{u^m \, m^m}{m!} \, \A(m,um)\,.
$$
\end{Lemma}

\begin{proof}
The conjugation,
\begin{equation}
e^{u \F_2} \, \al_{-m} \, e^{-u \F_2} = \, \cE_{-m} (u m) \,,\label{F2E}
\end{equation}
is easily calculated from the definitions since the operator
$e^{u \F_2}$ acts diagonally.

The operators $\cE$ satisfy the following basic commutation
relation:
\begin{equation}\label{commcE}
  \left[\cE_a(z),\cE_b(w)\right] =
\cs\left(\det  \left[
\begin{smallmatrix}
  a & z \\
b & w
\end{smallmatrix}\right]\right)
\,
\cE_{a+b}(z+w)\,.
\end{equation}
{}From \eqref{commcE}, we obtain
$$
[\al_1, \cE_{-m}(s)] = \cs(s) \ \cE_{-m+1}(s)
$$
and, therefore,
\begin{equation}
\label{ccc}
e^{\al_1} \, \cE_{-m}(s) \, e^{-\al_1}  =
\frac{\cs(s)^{m}}{m!}
\sum_{k\in \Z} \frac{\cs(s)^{k}}{(m+1)_k} \, \cE_k(s)
 \,.
\end{equation}
Applying \eqref{ccc} to \eqref{F2E} completes the proof.
\end{proof}

\subsubsection{}

Equation \eqref{lllst} and Lemma \ref{conjF} together yield a
concise formula for the evaluations of $\bH(z_1, \ldots, z_n,u)$ at the
positive integers
$z_i=\mu_i$:
\begin{equation}
\bH(\mu_1,\ldots,\mu_n,u) = u^{-n} \,
\lang \prod_{i=1}^n \A(\mu_i,u\mu_i) \rang \,.
\label{Hmat}
\end{equation}
However, we will require a stronger result.
We will prove that the right side of equation
\eqref{Hmat} is an analytic function of the variables $\mu_i$
and that the $n$-point function $\bH(z_1, \ldots, z_n,u)$
is a Laurent expansion of this analytic function.

\subsection{Convergence of matrix elements}

\subsubsection{}

If $a\ne 0,1,2,\dots$, the sum in \eqref{defA} is infinite
in both directions. Hence, for general values of  $\mu_i$,
the matrix
element on the right side of \eqref{Hmat} is not \`a priori
well-defined.
By expanding the definition of $\A(\mu_i, u\mu_i)$,
the right side
of \eqref{Hmat} is an $n$-fold series. We will prove the series
converges in a suitable domain of values of $\mu_i$.

Let $\Omega$ be the following domain in $\com^n$:
$$
\Omega=\left\{(z_1,\dots,z_n)\left|
|z_k| > \sum_{i=1}^{k-1} |z_i|,\,k=1,\dots,n \right. \right\}\,.
$$
The constant term of the operator ${\mathcal E}_0(uz_i)$ occurring in
the definition of $\A(z_i,uz_i)$ has a pole at $uz=0$.
For $u\neq 0$,
the coordinates $z_i$ are kept away in $\Omega$ from the poles
$uz_i=0$. 
We will prove the following convergence result.

\begin{Proposition}\label{Pc}
Let $K$ be a compact set,
$$
K\subset \Omega \cap \{z_i\ne -1,-2,\dots, i=1,\dots,n\}.
$$
For all partitions $\nu$ and $\lambda$, the series
\begin{equation}
\left( \A(z_1,uz_1) \cdots \A(z_n,u z_n) \, v_\nu, v_\la
\right)
\label{Am}
\end{equation}
converges absolutely and uniformly on $K$ for all sufficiently
small $u\neq 0$.
\end{Proposition}

\subsubsection{}
We will require three Lemmas for the proof of Proposition \ref{Pc}.

\begin{Lemma}\label{Lc1}
Let $\nu$ be a partition of $k$. For any integer $l$,
there exists at most $\max(k,l)$
partitions $\la$ of $l$ satisfying
$$
(\A(z,uz) \, v_\nu,v_\la) \ne 0 \,.
$$
\end{Lemma}

\begin{proof}
If $k=l$, then by the definition of $\A(z,uz)$,
there is exactly one such partition $\la$,
namely $\la=\nu$.

Next, consider the case $k>l$.
If the matrix element does not vanish,
then the operator $\cE_{k-l}$ in \eqref{defA} must
act on one of the factors of
$$
v_\nu = \ul{\nu_1-\tfrac12} \wedge   \ul{\nu_2-\tfrac32} \wedge
\ul{\nu_3-\tfrac52} \wedge \dots\,,
$$
and decrease the corresponding part of the partition $\nu$.
Since $\nu$ has at most $k$ parts, the above action can occur
in at most $k$ ways. The argument in the
$l>k$ case is similar.
\end{proof}

\begin{Lemma}\label{Lc2}
For any two partitions $\nu$ and $\la$ satisfying
$|\nu|\ne |\la|$, we have
$$
\left|(\cE_{|\nu|-|\la|}(uz) \, v_\nu,v_\la)\right| \le
\exp\left(\frac{|\nu|+|\la|}2 \, |uz|\right) \,.
$$
If $\nu=\la$, then
$$
\left|\left(\cE_0(uz)\, v_\nu,v_\nu \right)-
\frac1{\cs(uz)}
\right| \le |\nu| \, \exp(|\nu| |uz|) \,.
$$
\end{Lemma}
\begin{proof}
The Lemma is obtained from the definition of $\mathcal{E}_{|\nu|-|\lambda|}(uz)$.
\end{proof}

\begin{Lemma}
For all fixed $k_0,k_n\in \Z$, the series
\begin{equation}
  \label{ec2}
 \sum_{k_1,\dots,k_{n-1} \ge 0} \, \prod_{i=1}^n \,
\frac{z_i^{k_i-k_{i-1}}}{(d_i)_{k_i-k_{i-1}}}
\end{equation}
converges absolutely and uniformly on compact subsets of $\Omega$
for all values of the parameters $d_i\ne 0,-1,-2,\dots$.
\end{Lemma}

By differentiating with respect to the
variables $z_i$, we can insert in \eqref{ec2} any polynomial
weight in the summation variables $k_i$.

\begin{proof}
Consider the factor obtain by summation with respect to $k_1$:
\begin{equation} \label{ec3}
 \sum_{k_1\ge 0} \frac{(z_1/z_2)^{k_1}}{(d_1)_{k_1-k_0} (d_2)_{k_2-k_1}} \,.
 \end{equation}
The above series converges absolutely and uniformly
on compact sets since $|z_1/z_2|<1$ on the domain $\Omega$.
We require a bound on \eqref{ec3} considered as a function of
the parameter $k_2$.

The series \eqref{ec3} is bounded by a high enough derivative
of the series
\begin{equation} \label{ec4}
 \sum_{k_1\ge 0}^{k_2} \frac{w^{k_1}}{k_1!\, (k_2-k_1)!} +
\sum_{k_1 > k_2} \frac{(k_1-k_2)!}{k_1!} \, w^{k_1} \,,
\quad
w=\left|\frac{z_1}{z_2}\right| \,.
\end{equation}
The first term of \eqref{ec4} can be obviously estimated by
$$
\frac1{k_2!} \, \left(\frac{|z_1|+|z_2|}{|z_2|}\right)^{k_2} \,,
$$
whereas the second term of \eqref{ec4}
can be estimated by
$$
\frac1{k_2!} \, \frac{|z_1/z_2|^{k_2+1}}{1-|z_1|/|z_2|} \,.
$$
Therefore, the sum over both $k_1$ and $k_2$ behaves like the series
$$
\sum_{k_2\ge 0} \frac1{k_2! \,(d_3)_{k_3-k_2}}
\left(\frac{|z_1|+|z_2|}{|z_3|}\right)^{k_2} \,,
$$
which is a sum of the form \eqref{ec3}. Again, the series converges
absolutely and uniformly on compact sets since
$|z_1|+|z_2| < |z_3|$.

The Lemma  is proved by iterating the above argument.
\end{proof}

\subsubsection{Proof of Proposition \ref{Pc}}

We first expand \eqref{Am} as a sum over all intermediate
vectors
$$
v_\nu=v_{\mu[0]}, v_{\mu[1]},
\dots, v_{\mu[n-1]}, v_{\mu[n]} = v_\la.
$$
Next, using Lemmas \ref{Lc1} and \ref{Lc2}, we will  bound
the summation over all intermediate partitions $\mu$ by
a summation over their sizes,
$$
k_i = \left|\mu[i]\right| \,, \quad i=0,\dots, n \,.
$$
The term $\max(k_i,k_{i+1})$ of  Lemma \ref{Lc1} can be
bounded by $k_i+k_{i+1}$ and, in any case, amounts to
an irrelevant polynomial weight.

We conclude the Proposition will be established if
the absolute
convergence for $z\in K$ and sufficiently small $u$
of the following series is proven:
\begin{equation}
  \label{ec1}
 \sum_{k_1,\dots,k_{n-1} \ge 0} \, \prod k_i^{m_i} \,
e^{(k_i+k_{i-1})|uz_i|/2} \, \frac{\cs(uz_i)^{k_i-k_{i-1}}}{(1+z_i)_{k_i-k_{i-1}}}
\,,
\end{equation}
where the parameters $m_i$ are fixed nonnegative integers.
Here, we neglect the prefactors $\cS(uz_i)^{z_i}$ of 
the operators $\A(z_i, uz_i)$ ---
the functions $\cS(uz_i)^{z_i}$ are 
analytic and single valued (for the principal branch)
on $K$ for sufficiently small $u$. Also, 
we neglect
the constant terms of $\A(z_i,uz_i)$ as they
do not affect convergence for $u\neq 0$.

The terms raised to the power $k_i$ in
\eqref{ec1} are
$$
\left( e^{(|uz_i| + |u z_{i+1}|)/2} \, \frac{\cs(uz_i)}{\cs(uz_{i+1})}
\right)^{k_i} \,, \quad i=1,\dots, n-1 \,.
$$
Since, for $u\to 0$, we have
$$
e^{(|uz_i| + |u z_{i+1}|)/2} \, \frac{\cs(uz_i)}{\cs(uz_{i+1})}
\to \frac{z_i}{z_{i+1}} \,,
$$
the convergence of the series \eqref{ec1} follows from the
convergence of the series \eqref{ec2} with values
$$
d_i = 1 + z_i \,, \quad i =1 , \dots,  n \,.
$$
\qed

\subsection{Series expansion of matrix elements}

\subsubsection{}

By Proposition \ref{Pc}, the vacuum matrix element
\begin{equation}
\lang \A(z_1,uz_1) \cdots \A(z_n,u z_n)
\rang
\label{Am0}
\end{equation}
is an analytic function of the variables $z_1,\ldots, z_n,u$
on a punctured open set of
$\Omega\times 0$ in $\Omega \times \com^*$. Therefore, we may expand \eqref{Am0}
in a
convergent  Laurent power series.

First, viewing $u$ as a parameter, we expand in Laurent series in the
variables $z_1, \ldots, z_n$ in the following manner. 
For any point
$(z_2,\dots,z_n)$ in the domain
$$
\Omega'=\left\{(z_2,\dots,z_n)\left|
|z_k| > \sum_{i=2}^{k-1} |z_i|,\,k=2,\dots,n \right. \right\}\,,
$$
the function \eqref{Am0} is analytic and single-valued for $z_1$ in a
sufficiently small punctured neighborhood of the origin. Hence, the function
can be expanded there in a convergent Laurent
series. Every coefficient of that Laurent expansion is
an analytic function on the domain $\Omega'$ and, by iterating
the same procedure, can be expanded completely into a Laurent power
series.
The coefficients of the Laurent expansion in the variables $z_1, \ldots, z_n$ 
may be expanded as
Laurent 
series in $u$. 

Alternatively, we may expand the
function \eqref{Am0} in the variable $u$ first. Then,
the coefficients of the expansion are analytic functions
on the domain $\Omega$.

Later, we will identify the Laurent series expansion  of \eqref{Am0}
with the series $u^n\, \bH(z_1,\ldots,z_n,u)$.

\subsubsection{}

For any ring $R$,
define the ring  $R((z))$ by
$$
R((z)) = \left\{
\left.\sum_{i\in \Z} r_i z^i\right| r_i\in R, \, r_n=0,\, n\ll 0\right\} \,.
$$
In other words, $R((z))$ consists of formal Laurent
series in $z$ with coefficients in $R$ and exponents bounded
from below.

\begin{Proposition}\label{Pfs}
We have
$$
\lang \A(z_1,uz_1) \cdots \A(z_n,u z_n)
\rang
 \in \Q[u^{\pm 1}] ((z_n)) \dots ((z_1))\,.
$$
\end{Proposition}

\begin{proof}
The result follows by induction on $n$ from the
following property of the operators $\A$:
$$
\left(\A(z,uz)-\frac1{uz}\right)^* \, v_\mu = O\left(z^{-|\mu|}\right) \,.
$$
Indeed, with the exception of the term $(uz)^{-1}$ which
appears in the constant term of operator $\mathcal{E}_0(uz)$,
terms contributing to the
coefficient
$$
\left[z^{-k}\right]\,
\left(\A(z,uz)-\frac1{uz}\right)^*
$$
lower the energy by at least $k$ and, since
there are no vectors of negative energy, annihilate $v_\mu$ if
$k>|\mu|$\,.
\end{proof}

\subsubsection{}

Let $\A_k$ be the coefficients of the expansion of
the operator $\A(z,uz)$ in powers of $z$:
\begin{equation}
  \label{defAk}
  \A(z,uz)= \sum_{k\in \Z} \A_k\, z^k \,.
\end{equation}
As observed in the proof of Proposition \ref{Pfs},
the operator $\A_k$ for $k\ne -1$ involves only terms
of energy $\ge - k$. The same is true for $\A_{-1}$ with
the exception of the constant term $-u^{-1}$.

In terms of the operators $\A_k$, the Laurent
series expansion of \eqref{Am0} can be written
as:
\begin{equation}
  \label{Ams}
 \lang \A(z_1,uz_1) \cdots \A(z_n,u z_n)
\rang =
\sum_{k_1,\dots,k_n}  \lang \A_{k_1} \cdots \A_{k_n}
\rang
\,  z_1^{k_1}\ldots z_n^{k_n}  \,.
\end{equation}
If $k_j < -\sum_{i<j} (k_i+1)$  for some $j$, then
the corresponding term vanishes by energy
considerations.

\subsection{Commutation relations and rationality}

\subsubsection{}

Consider the doubly infinite series:
$$
\de(z,-w) = \frac1w \sum_{n\in\Z} \left(-\frac zw\right)^n  \ \in \Q((z,w)).
$$
The above series is the difference between  the following two
expansions:
\begin{align}\label{zw1}
  \frac1{z+w} &= \frac1w - \frac{z}{w^2} + \frac{z^2}{w^3} -
  \dots \,, \quad |z|<|w| \,, \\ \label{zw2}
\frac1{z+w} &= \frac1z - \frac{w}{z^2} + \frac{w^2}{z^3} -
  \dots \,, \quad |z|>|w| \,.
\end{align}
The series $\de(z,-w)$ is a formal $\delta$-function at $z+w=0$,
in the sense that
$$
(z+w) \, \de(z,-w) = 0 \,.
$$

\subsubsection{}

The following important result will be established in
Section \ref{ffff}.

\begin{Theorem}\label{ThA}
 We have
 \begin{equation}
   \label{cA1}
   [\A(z,uz),\A(w,uw)]  = zw\, \de(z,-w) \,,
 \end{equation}
or equivalently,
\begin{equation}
  \label{cA2}
  \left[\A_k, \A_l\right] = (-1)^l \delta_{k+l-1} \,.
\end{equation}
\end{Theorem}

\begin{Corollary}\label{cThA}
The series
\begin{equation}
\prod_{i<j} (z_i+z_j) \, \lang \A(z_1,uz_1) \cdots \A(z_n,u z_n)
\rang \in  \Q[u^{\pm1}] ((z_n))\dots ((z_1))
\label{cle}
\end{equation}
is symmetric in  $z_1,\dots,z_n$ and, hence, is an element
of
$$
\prod z_i^{-1} \, \Q[u^{\pm1}] \, [[z_1,\dots,z_n]] \,.
$$
\end{Corollary}

\begin{proof}
Indeed, the exponents of $z_1$ in \eqref{cle} are bounded below
by $-1$.
\end{proof}

\subsubsection{}

We now deduce the following result from Theorem \ref{ThA}:

\begin{Proposition}\label{twww}
The coefficients,
\begin{equation}\label{Amu}
\left[u^m\right]\,
\langle \A(z_1,uz_1) \dots \A(z_n,u z_n)\rangle \,, \quad m\in \Z\,,
\end{equation}
of powers of $u$ in the expansion \eqref{Ams}
are symmetric rational functions in $z_1,\dots,z_n$,
with at most simple poles on
the divisors $z_i+z_j=0$ and $z_i=0$.
\end{Proposition}

\begin{proof}
By Corollary \ref{cThA}, it suffices to prove
the exponents of $z_n$ in the expansion of \eqref{Amu}
are bounded from above.

The equation,
\begin{equation} \label{re1}
  \lang
\cE_{k_1}(uz_1)\ldots \cE_{k_n}(uz_n)
\rang = \lang
\frac{\cE_{k_1}(uz_1)}{u^{k_1}}\ldots \frac{\cE_{k_n}(uz_n)}{u^{k_n}}
\rang,
\end{equation}
holds
since the vacuum expectation vanishes unless $\sum k_i = 0$.
The transformation $\cE_k \to u^{-k} \cE_k$
applied to the operator $\A(z,uz)$ acts as the
substitution
$$
\cs(uz)^k \mapsto \frac{\cs(uz)^k}{u^k} \,,
$$
which makes all terms regular and nonvanishing at $u=0$,
except for the simple pole in the constant term $\cs(uz)^{-1}$.

Since \eqref{re1} vanishes if $k_n >0$, the vacuum
expectation
\begin{equation*}
\langle \A(z_1,uz_1) \dots \A(z_n,u z_n)\rangle
\end{equation*}
 depends on $z_n$ only through
terms of the form
$$
\cS(uz_n)^{z_n}  \,, \quad
e^{auz_n} \,, \quad a\in \tfrac12\Z\,,
$$
as well as
\begin{multline*}
  z_n(z_n-1)\dots(z_n-k+1) \, \frac{u^k}{\cs(uz_n)^{k}} = \\
  \left(1-\frac{1}{z_n}\right) \cdots \left(1-\frac{k-1}{z_n}\right)
  \, \cS(uz_n)^{-k}\,, \quad k=1,2,\dots \,.
\end{multline*}
Because these terms are multiplied by a function of $u$
with a bounded order of pole at $u=0$,
the required boundedness of degree in $z_n$ for fixed powers of $u$
is now immediate.
\end{proof}

\subsection{Identification of $\bH(z,u)$}

\subsubsection{}

By definition \eqref{bHu}, $\bH(z_1,\dots,z_n,u)$ is a Laurent series
in $u$ with coefficients given by rational functions
of $z_1,\dots,z_n$ which have at most first order poles
at the divisors $z_i+z_j=0$ and $z_i=0$.

By Proposition \ref{twww}, the expansion \eqref{Ams}
has the exact same form.
We can now state the main result of the present
section.

\begin{Theorem} \label{hodfor} We have
\begin{equation}
\bH(z_1, \ldots, z_n,u) = u^{-n}
\lang \A(z_1,uz_1) \dots \A(z_n,u z_n) \rang \,. \label{HA}
\end{equation}
\end{Theorem}

\begin{proof}
By Proposition \ref{Pc}, the coefficients \eqref{Amu} are
analytic functions on the domain $\Omega$. Moreover, by
Proposition \ref{twww}, these functions are rational.
By \eqref{Hmat}, for positive integer values of $z_i$
in the domain $\Omega$, these functions take the same values
as the corresponding coefficients of $\bH$. Since
positive integer values of $z_i$ inside $\Omega$ form a Zariski
dense set, the Theorem follows.
\end{proof}

\subsubsection{}

As an illustration of Theorem \ref{hodfor}, we obtain
the following result.

\begin{Proposition}
  The connected $2$-point generating function
$\bH^\circ(z_1,z_2,u)$ for Hodge integrals
is given by
\begin{multline}
  \label{H2c}
  \bH^\circ(z_1,z_2,u)
= \frac{\cS(uz_1)^{z_1}\, \cS(uz_2)^{z_2}}{\cs(u(z_1+z_2))}\times \\
\left[\Ff{-z_2}{1+z_1}{\frac{1-e^{uz_1}}{1-e^{-uz_2}}} -
 \Ff{-z_2}{1+z_1}{\frac{1-e^{-uz_1}}{1-e^{uz_2}}}  \right] \,,
\end{multline}
where ${}_2F_1$ the Gauss hypergeometric function  \eqref{fS}.
\end{Proposition}

\begin{proof}
We first calculate:
\begin{multline}\label{EEv}
\lang \cE_{k_1}(uz_1) \, \cE_{k_2}(uz_2) \rang -
 \lang \cE_{k_1}(uz_1) \rang \, \lang\cE_{k_2}(uz_2) \rang
= \\
\begin{cases}
\dfrac{\cs(k_1 u (z_1+z_2))}{\cs(u (z_1 + z_2))}\,,
& 0 < k_1 = - k_2\,, \\
0\,, & \textup{otherwise.}
\end{cases}
\end{multline}
The nonzero term in \eqref{EEv} arises from the constant term
in the commutator $\left[\cE_{k_1}(uz_1), \cE_{-k_1}(uz_2)\right]$.
Then, by formula \eqref{HA}, we obtain
\begin{multline*}
  \bH^\circ(z_1,z_2,u)
= \\
\frac{\cS(uz_1)^{z_1}\, \cS(uz_2)^{z_2}}{\cs(u(z_1+z_2))}\,
\sum_{k>0} \cs(k u (z_1+z_2)) \,
\frac{\cs(uz_1)^k \, \cs(uz_2)^{-k}}{(1+z_1)_k \, (1+z_2)_{-k}}\,,
\end{multline*}
which is equivalent to \eqref{H2c} \,.
\end{proof}

The symmetry in $z_1$ and $z_2$ is not at all
obvious from formula \eqref{H2c}.



\section{Operator formula for Gromov-Witten invariants}
\label{ththth}

\subsection{Localization revisited}

\subsubsection{}

Propositions \ref{Loc} and \ref{hodfor} together yield
the following localization formula in terms of vacuum expectations:
\begin{multline}\label{Gd}
  \bG_d(z_1,\dots,z_n,w_1,\dots,w_m,u)= \\
\sum_{|\mu|=d} \, \frac{1} {\zz(\mu)} \,
\bJ(z,\mu,u,t) \, \bJ(w,\mu,u,-t) \,,
\end{multline}
where the function $\bJ(z,\mu,u,t)$ is defined by:
\begin{multline}\label{defJ}
  \bJ(z_1,\dots,z_n,\mu_1, \ldots, \mu_\ell,u,t)= \\
t^{-d} u^{-n} \left(\prod
    \frac{\mu_i^{\mu_i}}{\mu_i !}\right) \,
\lang \prod  \A(t z_i,u z_i)
\prod \A(\mu_i,\tfrac ut \mu_i) \rang = \\
u^{-d-n}  \,
\lang \prod  \A(t z_i,u z_i)
\, e^{\al_1} \, e^{\frac ut \F_2} \, \prod \al_{-\mu_i} \rang
\, .
\end{multline}

\subsubsection{}

For each partition $\mu$, define the
vector $\chi_\mu \in \LV$ by:
$$
\chi_\mu =  \prod_{i=1}^{\ell(\mu)} \al_{-\mu_i} \, \vac \,.
$$
The expansion of $\chi_\mu$
in the standard basis $v_\nu$ is given
by the values of the symmetric group characters $\chi^\nu$ on the
conjugacy class determined by $\mu$:
$$
\chi_\mu =  \sum_{|\nu|=|\mu|}\chi^\nu_\mu \, v_\nu \,.
$$
From the
commutation relations
\begin{equation}
  \label{commal}
  \left[\al_k,\al_l\right] = k \, \delta_{k+l} \,,
\end{equation}
or from
the orthogonality relation for characters, we find
$$
(\chi_\mu,\chi_\nu)= \zz(\mu)\,  {\delta_{\mu,\nu}}.
$$

Let $\Pv$ denote the orthogonal projection onto the vector
$\vac$. Since the vectors $\{\chi_\mu\}_{|\mu|=d}$
span the eigenspace of $H$ with eigenvalue $d$, the operator
$$
P_d = \sum_{|\mu|=d} \,
\frac{1} {\zz(\mu)} \,
\prod \al_{-\mu_i} \, \Pv \, \prod \al_{\mu_i}
$$
is the orthogonal projection onto the $d$-eigenspace of $H$.

\subsubsection{}

Using definition \eqref{defJ} and
the projection $\Pv$, we can write
\begin{multline*}
  u^{2d+n+m}
\bJ(z,\mu,u,t) \, \bJ(w,\mu,u,-t) = \\
\lang \prod  \A(t z_i,u z_i)
\, e^{\al_1} \, e^{\frac ut \F_2} \, \prod \al_{-\mu_i} \,
\Pv\, \times  \right. \\
\left.
\prod \al_{\mu_i} \, \, e^{-\frac ut \F_2} \, e^{\al_{-1}} \,
\prod  \A(-t w_j,u w_j)^*
\rang
\,.
\end{multline*}
Since  $\F_2$ commutes with $H$,  $\F_2$ also
commutes with $P_d$. Therefore,
\begin{equation}\label{Pd}
\sum_{|\mu|=d}
\frac{1} {\zz(\mu)} \,
e^{\frac ut \F_2}
\prod \al_{-\mu_i} \, \Pv \, \prod \al_{\mu_i}   e^{-\frac ut \F_2} =
P_d
\end{equation}
After summing \eqref{Gd} using \eqref{Pd},
 we find:
\begin{multline}\label{Gdop}
  \bG_d(z_1,\dots,z_n,w_1,\dots,w_m,u)= \\
u^{-2d-n-m} \lang \prod  \A(t z_i,u z_i)
\, e^{\al_1} \, P_d \, e^{\al_{-1}} \,
\prod  \A(-t w_j,u w_j)^*
\rang \,.
\end{multline}

\subsubsection{}

Define the $n+m$-point function $\bG(z,w,u)$ of equivariant Gromov-Witten
invariants of all degrees by:
$$
\bG(z,w,u)= \sum_{d \ge 0} q^{d} \, \bG_d (z,w,u)\,.
$$
Since $H=\sum_{d} d \, P_d$, we find:
\begin{multline}\label{Gm}
  \bG(z_1,\dots,z_n,w_1,\dots,w_m,u)= \\
 u^{-n-m}
\lang \prod  \A(t z_i,u z_i)
\, e^{\al_1} \, \left(\frac{q}{u^2}\right)^H \, e^{\al_{-1}} \,
\prod  \A(-t w_j,u w_j)^*
\rang \,.
\end{multline}
Introduce the following operators:
\begin{align} \label{bA}
   \bA(z) &= \frac1u \, \A(t z,u z), \\
\bA^\star(w) &= \frac1u \, \A(-t w,u w)^* \notag \,.
\end{align}
Recall, by definition,
\begin{equation}
  \label{dbA}
  \bA(z) = u^{-1} \, \cS(uz)^{tz} \,
\sum_{k\in \Z} \frac{\cs(uz)^{k}}{(1+tz)_k} \, \cE_k(uz)
 \,.
\end{equation}
We obtain the following result by substituting the operators
$\bA(z)$, $\bA^\star(w)$ in equation \eqref{Gm}.

\begin{Theorem}\label{twoo}
 The function $\bG(z,w,u)$ is the following vacuum expectation:
\begin{multline}\label{bG}
  \bG(z_1,\dots,z_n,w_1,\dots,w_m,u)= \\
\lang \prod  \bA(z_i) \,
\, e^{\al_1} \, \left(\frac{q}{u^2}\right)^H \, e^{\al_{-1}} \,
\prod  \bA^\star(w_j)
\rang
 \,.
\end{multline}
\end{Theorem}

In particular, for the $0$-point function, Theorem \ref{twoo} yields the
following correct evaluation:
$$
\bG() = \lang
\, e^{\al_1} \, \left(\frac{q}{u^2}\right)^H \, e^{\al_{-1}}
\rang = e^{q/u^2} \,.
$$

\subsection{The $\tau$-function}

\subsubsection{}

By definition, $\bG(z,w,u)$ includes unstable contributions obtained
from \eqref{gunst}. We will now introduce the $\tau$-function: a
generating
function for the true equivariant Gromov-Witten invariants of $\proj^1$.
The $\tau$-function does not include unstable
contributions.
In Theorems \ref{TT1} and \ref{TT2}, we will show the $\tau$-function of
the equivariant theory of $\proj^1$ is
a $\tau$-function of an integrable hierarchy, namely,
the 2--Toda hierarchy of Ueno and Takasaki.

\subsubsection{}

Let $\bA_k$ denote the coefficient of $z^{k+1}$ in the
expansion of $\bA$:
$$
\bA_k = [z^{k+1}] \, \bA \,, \quad
\bA^\star_k = [z^{k+1}] \, \bA^\star \,, \quad
k\in \Z\,.
$$
Then, by Theorem \ref{twoo},
\begin{multline}\label{pedd}
  \sum_{g\in \Z} \sum_{d\geq 0} u^{2g-2} q^d
\lang \prod \tau_{k_i}(\bz) \prod \tau_{l_j}(\bi) \rang^\bullet_{g,d} =\\
\lang \prod  \bA_{k_i} \,
\, e^{\al_1} \, \left(\frac{q}{u^2}\right)^H \, e^{\al_{-1}} \,
\prod  \bA^\star_{l_j}
\rang
 \,,
\end{multline}
where, the left side consists of  true equivariant Gromov-Witten
invariant (with no unstable contributions). The unstable contributions
\eqref{gunst} produce terms of degrees at most 0 in their variables and, therefore,
do not contribute to \eqref{pedd}.

\subsubsection{}

Let the variable sets $x_i, x^\star_i$ correspond to the descendents
$\tau_i(\bz), \tau_i(\bi)$ respectively.
Define the equivariant  $\tau$-function by:
$$
\tau(x,x^\star,u) =
\sum_{g\in \Z} \sum_{d\geq 0} u^{2g-2} \, q^d \,
\lang \exp\left(\sum_{i\ge 0}
x_i \, \tau_i(\bz)+x^\star_i \, \tau_i(\bi)
\right) \rang^\bullet_{g,d} \,.
$$

\begin{Theorem} \label{vvv} The equivariant $\tau$-function is a vacuum expectation in $\LV$:
  \begin{equation}
    \label{Z}
  \tau(x,x^\star,u) =   \lang e^{\sum x_i \bA_i} \,
\, e^{\al_1} \, \left(\frac{q}{u^2}\right)^H \, e^{\al_{-1}} \,
e^{\sum x^\star_i \bA^\star_i}
\rang \,.
  \end{equation}
\end{Theorem}
\begin{proof} The formula is a restatement of (\ref{pedd}).
\end{proof}

\subsection{The GW/H correspondence}

The generating function for the absolute stationary non-equivariant
Gromov-Witten theory of $\proj^1$ is obtained from
the generating function \eqref{bG} by taking
$$
m=0\,, \quad t=0\,, \quad u=1 \,.
$$
The operator formula \eqref{bG} then specializes to
$$
\bG(z,\emptyset,1)\big|_{t=0} = \lang \prod  \A(0,z_i) \,
\, e^{\al_1} \, q^H \, e^{\al_{-1}} \rang
 \,.
$$
We have
\begin{align}
  \A(0,z) & = \sum_{k \ge 0} \frac{\cs(z)^{k}}{k!} \, \cE_k(z)  \notag\\
& = e^{\al_1} \, \cE_0(z) \, e^{-\al_1} \,,
\end{align}
where the second equality follows from \eqref{ccc}.
We obtain the following result.

\begin{Proposition} The $n$-point function of
absolute stationary non-equivariant Gromov-Witten invariants of $\proj^1$
is given by:
\begin{equation}\label{gwh1}
\bG(z,\emptyset,1)\big|_{t=0} =
\lang e^{\al_1} \, q^H \, \prod \cE_0(z_i) \,  e^{\al_{-1}} \rang
 \,.
\end{equation}
\end{Proposition}

Extracting the coefficient of $q^d$ in \eqref{gwh1}, we obtain
the following equivalent formula:
\begin{equation}\label{gwh2}
\bG_d(z,\emptyset,1)\big|_{t=0} =
\frac{1}{(d!)^2}
\lang \al_1^d \, \prod \cE_0(z_i) \,  \al_{-1}^d \rang
 \,.
\end{equation}
This is precisely the special case of the GW/H correspondence
\cite{OP2} required for the proof of the general GW/H
correspondence given in \cite{OP2}.



\section{The 2--Toda hierarchy}
\label{fififi}

\subsection{Preliminaries on the 2--Toda hierarchy}

\subsubsection{}

Let $M$ be an element of the group $GL(\infty)$ acting
in the $GL(\infty)$-module $\LV$. The matrix elements
of the operator $M$,
$$
\left( M v, w\right)\,, \quad v,w\in \LV\,,
$$
can be viewed as, suitably regularized,
 $\frac\infty 2 \times \frac\infty 2$-minors
of the matrix $M$. In particular, the matrix elements satisfy
quadratic Pl\"ucker relations.

A concise way to write all the Pl\"ucker relations
is the following, see for example \cite{Kac,MJD}. Introduce the
following operator on $\LV\otimes \LV$:
$$
\Omega=\sum_{k\in  \Z+ \frac12} \psi_k \otimes \psi^*_k \,.
$$
The operator $\Omega$
 operator can be defined $GL(\infty)$-invariantly
by taking, instead of $\{\psi_k\}$ and $\{\psi_k^*\}$,
any linear basis of the space $V$ of creation
operators and the corresponding dual basis
of the space of annihilation operators.
The $GL(\infty)$-invariance implies
\begin{equation}
  \label{Pl1}
 \left[M\otimes M, \Omega\right] = 0
\end{equation}
for any operator $M$ in the closure of the image of
$GL(\infty)$ in the endomorphisms of $\LV$.

Concretely, for any $v,v',w,w'\in \LV$, we obtain the
following quadratic relation between the matrix
coefficients of $M$
\begin{equation}
  \label{Pl2}
  \left(\left[M\otimes M, \Omega \right]  \,v\otimes v',
w\otimes w'\right) = 0\,.
\end{equation}

\subsubsection{}

For example, consider the following vectors in \eqref{Pl2}:
\begin{align*}
  v &= \vac = \ul{-\tfrac12} \wedge \ul{-\tfrac32} \wedge \ul{-\tfrac52}
\wedge \ul{-\tfrac72} \wedge \dots  \,, \\
v' &= v_\square = \ul{\tfrac12} \wedge \ul{-\tfrac32} \wedge \ul{-\tfrac52}
\wedge \ul{-\tfrac72} \wedge \dots\,, \\
w &= v_1 = \ul{\tfrac12} \wedge
\ul{-\tfrac12} \wedge \ul{-\tfrac32} \wedge \ul{-\tfrac52}
\wedge \dots  \,, \\
w' &= v_{-1} = \ul{-\tfrac32} \wedge \ul{-\tfrac52}
\wedge \ul{-\tfrac72} \wedge \ul{-\tfrac92}\wedge \dots  \,,
\end{align*}
where $\vac,v_1,v_{-1}$ are the vacua in subspaces of charge
$0$, $1$, and $-1$, respectively, and $v_\square$ is the
unique charge 0 vector of energy 1, corresponding to the partition
$\lambda=(1)$.

We find from the definitions,
\begin{align*}
 \Omega \, \vac\otimes v_\square  &= v_1 \otimes v_{-1} \,, \\
\Omega^* \, v_1 \otimes v_{-1} &=  \vac\otimes v_\square -
 v_\square \otimes \vac \,.
\end{align*}
Hence, \eqref{Pl2} yields
the following identity:
\begin{multline}
  \label{Pl3}
  \left(M \, v_1, v_1\right) \, \left(M \, v_{-1}, v_{-1} \right) =  \\
\left(M \, \vac, \vac\right) \, \left(M \, v_\square, v_\square \right) -
\left(M \, \vac, v_\square \right) \,
\left(M \, v_\square, \vac \right) \,.
\end{multline}
The above identity, which remains valid for matrices of finite
size, is often associated with Lewis Carroll \cite{LC},
but was first established by P.~Desnanot in 1819 (see \cite{muir}).

Another way to write
identity \eqref{Pl3} is the following:
\begin{equation}
  \label{Pl4}
\lang T^{-1} M T \rang \, \lang T M T^{-1} \rang
= \lang M \rang \, \lang \al_{1} \, M \, \al_{-1} \rang -
\lang \al_{1} \, M  \rang \,
\lang M \, \al_{-1} \rang \,
\,,
\end{equation}
where $T$ is the translation operator on the
in infinite wedge space
$$
T \cdot \bigwedge \ul{s_i}  = \bigwedge \ul{s_i+1} \,.
$$

\subsubsection{} \label{taucon}

Using the \emph{vertex operators}
$$
\Gamma_\pm(t) = \exp\left(\sum_{k>0} t_k \, \frac{\al_{\pm k}}{k}
\right)\,,
$$
we define a  sequence
of \emph{$\tau$-functions} corresponding to
the operator $M$,
$$
\tau^M_n(t,s) = \lang T^{-n}\, \widehat{M} \,T^n\rang\,,
\quad \widehat{M} = \Gamma_+(t)  \, M \, \Gamma_-(s)\,, \quad
n\in \Z \,.
$$
The derivatives  of $\tau^M_n$ with respect to the
variables $t$ and $s$ are nothing but matrix elements of
the matrix $\widehat{M}\in GL(\infty)$. Hence, the
functions $\tau^M_n$ satisfy a collection of
bilinear partial differential equations. This collection
is known as the 2--Toda hierarchy of Ueno and Takasaki, see \cite{UT}
and also, for example, the Appendix to \cite{iw} for a brief
exposition.

In particular, the lowest equation of the hierarchy is a
restatement of the equation \eqref{Pl4}:
\begin{equation}
  \label{Toda1}
  \tau_n \,\frac{\partial^2}{\partial t_1 \partial s_1} \, \tau_n
- \frac{\partial}{\partial s_1} \, \tau_n \,
\frac{\partial}{\partial t_1} \, \tau_n = \tau_{n+1} \, \tau_{n-1} \,, \quad n\in \Z.
\end{equation}
We may rewrite \eqref{Toda1}  as:
\begin{equation}
  \label{Toda2}
 \frac{\partial^2}{\partial t_1 \partial s_1} \, \log\tau_n =
\frac{\tau_{n+1} \, \tau_{n-1}}{\tau_n^2} \,.
\end{equation}

\subsection{String and divisor equations}

\subsubsection{}

The equivariant divisor equations describes the
effects of insertions of $\tau_0(\bz)$ and $\tau_0(\bi)$.
In terms
of the disconnected $(n+m)$-point generating
function $\bG_d(z_1,\dots,z_n,w_1,\dots,w_m,u)$,
the divisor equation for $\tau_0(\bz)$
insertion takes the
following form.

\begin{Proposition}
  We have
  \begin{multline}
    \label{divis}
    \left[z_0^1\right] \, \bG_d(z_0,z_1,\dots,z_n,w_1,\dots,w_m,u) = \\
\left(d-\frac1{24}+t\sum_{i=1}^n z_i \right) \,
\bG_d(z_1,\dots,z_n,w_1,\dots,w_m,u)\,.
  \end{multline}
\end{Proposition}

Recall, by construction, the function $\bG_d$ includes
contributions from unstable moduli spaces. Therefore, the usual
geometric proof of the divisor equation requires a modification.
Instead, we will prove the formula \eqref{divis} using the
operator formalism.

The presence of the
disconnected and unstable
contributions in $\bG_d$ actually simplifies the form
of the divisor equation ---  special handling
of the exceptional cases is no longer required.

\begin{proof}

Equation \eqref{bG} states:
$$
\bG_d(z,w,u)= u^{-2d}
\lang \prod \bA(z_i) \, e^{\al_1} \, P_d \, e^{\al_{-1}} \,
\prod \bA(w_i)^\star \rang \,,
$$
and hence
$$
\left[z_0^1\right]\, \bG_d(z_0,z_1,\dots, z_n,w,u) =
\lang \bA_0 \, \prod \bA(z_i) \, e^{\al_1} \, P_d \, e^{\al_{-1}} \,
\prod \bA(w_i)^\star \rang \,.
$$

The operator $\bA_0$ has the following form
\begin{equation}
  \label{fA0}
  \bA_0 = \al_1 - \frac1{24} + \dots \,,
\end{equation}
where the dots stand for terms for which the adjoint
annihilates the vacuum.  Since the energy operator
$H$ also annihilates the vacuum, we can write:
\begin{multline}\label{al+H}
  \left[z_0^1\right]\, \bG_d(z_0,z_1,\dots,z_n,w,u) = \\ \lang
  \left(-\tfrac1{24}+\al_1+H\right) \, \prod \bA(z_i) \, e^{\al_1} \,
  P_d \, e^{\al_{-1}} \, \prod \bA(w_i)^\star \rang \,.
\end{multline}
From  definition \eqref{dbA}, we find:
\begin{equation}\label{divA}
  \left[\al_1+H,\bA(z)\right] = t z \, \bA(z)  \,.
\end{equation}
Also, we have $[H,\al_1]=-\al_1$ and $H\, P_d = d \, P_d$. Therefore,
$$
(\al_1+H) \, e^{\al_1} \, P_d  = e^{\al_1} \, H \, P_d = d \,
e^{\al_1} \, P_d \,.
$$
Hence, commuting the operator $\al_1+H$ in \eqref{al+H} to the middle,
we obtain formula \eqref{divis}.
\end{proof}

\subsubsection{}
The string equation describes the effect of the insertion of
$\tau_0(1)$, where $1$ is the identity class in the equivariant
cohomology of $\proj^1$. Since
$$
1 = \frac{\bz-\bi}t
$$
in the localized  equivariant
cohomology of $\proj^1$, the string equation is a linear combination
of the divisor equations associated to two torus fixed points.
The effect of an arbitrary number of the $\tau_0(1)$-insertions
can be conveniently described in the following form.

\begin{Proposition} We have
  \begin{multline}\label{string}
    \lang e^{\tau_0(1)} \, \prod \tau_{k_i}(\bz) \,
\prod \tau_{l_i} (\bi) \rang^\bullet_{g,d} = \\
\left[ \prod z_i^{k_i+1} \, \prod w_i^{l_i+1} \right] \,
e^{\sum z_i + \sum w_i} \, \bG_{g,d}(z,w,u) \,.
  \end{multline}
\end{Proposition}

\subsection{The 2--Toda equation}

\subsubsection{}

Let $\MMM$ be the matrix appearing in \eqref{Z},
\begin{equation}\label{dMMM}
\MMM = e^{\sum x_i \bA_i} \,
\, e^{\al_1} \, \left(\frac{q}{u^2}\right)^H \, e^{\al_{-1}} \,
e^{\sum x^\star_i \bA^\star_i} \,.
\end{equation}
In Section \ref{22TT}, we will see that for a suitable
matrix $M$, one can 
conjugate $\MMM$ to the canonical form $\Gamma_+(t) \, M \, \Gamma_-(s)$  
required of the 2-Toda hierarchy.  Here, the time
variables $\{t_i\}$ and $\{s_i\}$ are related to the variables
$\{x_i\}$ and $\{x^\star_i\}$ by an explicit linear transformation. 

The 2--Toda equation, the lowest equation of the Ueno-Takasaki
hierarchy, is then
a consequence of the results Section \ref{22TT}.
However, a direct derivation of the 2--Toda equation, without
the full hierarchy, is presented here first.

\subsubsection{}

{}From \eqref{fA0} we obtain
$$
\frac{\partial}{\partial x_0} \, \tau(x,x^\star,u) =
\lang (\al_1-\tfrac1{24}) \, \MMM \rang \,,
$$
and, similarly,
$$
\frac{\partial}{\partial x^\star_0} \, \tau(x,x^\star,u) =
\lang \MMM \, (\al_{-1}-\tfrac1{24})  \rang \,.
$$
We therefore find
\begin{multline}\label{TMT}
  \tau \, \frac{\partial^2}{\partial x_0 \, \partial x^\star_0} \,\tau -
\frac{\partial}{\partial x_0} \,\tau \,
\frac{\partial}{\partial x^\star_0} \,\tau  = \\
 \lang \MMM \rang \, \lang \al_{1} \, \MMM \, \al_{-1} \rang -
\lang \al_{1} \, \MMM  \rang \,
\lang \MMM \, \al_{-1} \rang  = \\
\lang T^{-1} \, \MMM \, T \rang \, \lang T \, \MMM \, T^{-1} \rang \,,
\end{multline}
where the second equality follows from \eqref{Pl4}\,.

\subsubsection{}

We will now study the conjugation of $\MMM$ by the translation
operator $T$.  
The result combined
with \eqref{TMT} will yield the 2--Toda equation.

We first examine the $T$ conjugation of the constituent 
operators of $\MMM$.
The conjugation of the operators
$\bA_k$ is best summarized by the equation
\begin{equation}
T^{-1} \, \bA(z) \, T = e^{uz}  \, \bA(z) \,,
\label{Ta2}
\end{equation}
which follows directly from definitions.
The conjugation equations for $\bA_k^\star$ are identical.

Since $T$ commutes with $\al_{\pm1}$, the only other
conjugation we require is:
\begin{equation}
T^{-n} \, H \, T^{n} = H + n  C + \frac{n^2}2 \,,
\label{Ta1}
\end{equation}
where $C$ is the charge operator
(see
Section 2.2.3 of \cite{OP2}).
Since
$C$ commutes with the 
remaining  operators $\bA_k, \bA_k^\star, \al_{\pm 1}$ 
and annihilates
the vacuum,  we may ignore $C$.

We now observe
the evolution of the
operators $\bA_k, \bA_k^\star$ under the string equation in
\eqref{string} has exactly same form as \eqref{Ta2}. Introduce,
the following differential operator
$$
\partial = \frac1{t} \, \left(\frac{\partial}{\partial x_0} -
\frac{\partial}{\partial x^\star_0}\right) \,,
$$
the action of which on $\tau$ corresponds to the
insertion of $\tau_0(1)$.

Combining \eqref{Ta2}, \eqref{Ta1}, and \eqref{string}, we obtain
\begin{equation}
  \label{TM}
  \lang T^{-n} \, \MMM \, T^n \rang =
\frac{q^{n^2/2}}{u^{n^2}} \, e^{nu\partial} \tau \,,
\end{equation}
and therefore,
\begin{equation*}
\lang T^{-1} \, \MMM \, T \rang \, \lang T \, \MMM \, T^{-1} \rang =
\frac{q}{u^2} \, e^{u\partial} \tau \, \, e^{-u\partial} \tau \,.
\end{equation*}
Thus, we have established the following version of the 2--Toda
equation for the function $\tau(x,x^\star,u)$.

\begin{Theorem}\label{TT1} The function $\tau(x,x^\star,u)$
satisfies the following form of the 2--Toda equation:
\begin{equation}
\frac{\partial^2}{\partial x_0\, \partial x^\star_0} \,
\log \tau = \frac{q}{u^2} \,
\frac{e^{u\partial} \tau \, \, e^{-u\partial} \tau}{\tau^2} \,.
\label{Toda}
\end{equation}
\end{Theorem}

Since the degree variable $q$ appears as a factor on the
right side of \eqref{Toda}, the equation \eqref{Toda}
determines all positive degree Gromov-Witten invariants of $\proj^1$
from the degree $0$ invariants.

\subsection{The 2--Toda hierarchy}\label{22TT}

\subsubsection{}\label{fullT}

Our goal now is to prove that there exists a upper
unitriangular matrix $W$ such that
\begin{equation}
W^{-1} \, \exp\left(\sum x_i \, \bA_i \right)  \, W = \Gamma_+(t)\,,
\label{WGa}
\end{equation}
where the time variables $\{t_i\}$ are obtained from the 
variables $\{x_i\}$ by certain explicit linear transformation
which will be described below. 

Once \eqref{WGa} is established, one deduces the 
2-Toda hierarchy for the $\tau$-function \eqref{Z}
as follows. First, taking the adjoint of the equation 
\eqref{WGa} and reversing the sign of the 
equivariant parameter $t$, we obtain
\begin{equation}
W^\star \, \exp\left(\sum x^\star_i \, \bA^\star_i \right)  \, 
\left(W^\star\right)^{-1} = \Gamma_-(s)\,
\label{WGas}
\end{equation}
where 
$$
W^\star=W^*\big|_{t\mapsto -t} \,. 
$$
The linear transformation 
$$
\{x^\star_i\}\mapsto\{s_i\}
$$ 
is obtained from 
the linear transformation $\{x_i\}\mapsto\{t_i\}$ 
by reversing the sign of the equivariant parameter $t$. 

Together, the equations \eqref{WGa} and \eqref{WGas}, give
the following formula for the matrix \eqref{dMMM}
\begin{equation}
  \MMM = W \, \Gamma_+(t) \, M \, \Gamma_-(s) \, W^\star \,,
\end{equation}
where
$$
M = W^{-1} \, 
e^{\al_1} \, \left(\frac{q}{u^2}\right)^H \, e^{\al_{-1}} \,
\left(W^\star\right)^{-1} \,.
$$
The unitriangularity of $W$ implies
\begin{equation*}
W^* \, \vac = W^\star \,\vac = \vac\,,
\end{equation*}
and, more generally,
\begin{equation*}
W^* \, T^n \, \vac = W^\star \, T^n \, \vac =  T^n \, \vac\,, \quad n\in \Z\,.
\end{equation*}
Therefore, we obtain 
\begin{align}
 \frac{q^{n^2/2}}{u^{n^2}} \, e^{nu\partial} \tau & =  
\lang T^{-n} \, \MMM \, T^n \rang \notag \\
& =
\lang T^{-n} \,\Gamma_+(t) \, M \, \Gamma_-(s)  \, T^n \rang \,,
\label{tautransf}
\end{align}
where the first equation is copied from \eqref{TM}. It then 
follows that the sequence \eqref{tautransf}
is a sequence of $\tau$-functions for
the full 2--Toda hierarchy of Ueno and Takasaki.

\subsubsection{} \label{stend}

We now proceed with the realization of the above plan. 

We will now view the operators $\bA_k$ as matrices
in the
associative algebra $\End(\infty)$. All
multiplication operations in Sections \ref{stend} -- \ref{endend}
should
be interpreted as multiplication in $\End(\infty)$, and {\em not} in
$\End(\LV)$.

For $k\geq 0$,
the matrices $\bA_k$ commute by Theorem \ref{ThA} and have
the form
\begin{equation} \label{fffe}
\bA_k = \frac{u^k}{(k+1)!} \, \al_{k+1} + \dots \,, 
\end{equation}
where the dots stand for term of energy larger than $-k-1$.

Since the matrix $\bA_0$ has form \eqref{fffe}, 
there exists an upper unitriangular matrix $W\in GL(\infty)$
conjugating $\bA_0$ to $\al_1$:
$$
W^{-1} \, \bA_0 \, W= \al_1\,.
$$
We call the matrix $W$ the \emph{dressing operator}.
The explicit form of $W$ is rather complicated,
unique only up to left multiplication
by a element of the centralizer of $\al_1$,
 and
will not be required.

However, the
dressed matrices
$$
\bAt_k = W^{-1}\, \bA_k \, W,  \quad k\geq 0,
$$
are uniquely defined and can be identified explicitly.

Because the matrices $\bAt_k$ commute with the
matrix $\bAt_0=\al_1$, the matrices have the following form:
\begin{equation}
  \label{WAWc}
  \bAt_k =
\sum_{l\le k+1} c_{k,l}(u,t) \, \al_l \,, \quad k=0,1,\dots\,,
\end{equation}
where
$$
c_{k,k+1} = \frac{u^k}{(k+1)!} \,.
$$
The other coefficients of the expansion
are determined by the following result.

\begin{Theorem}\label{ThW} The dressed operators $\bAt_k$
are determined by a generating function
identity:
\begin{equation}
\label{WAW}
\sum_{k\ge 0}  z^{k+1}\, \bAt_k =
\sum_{n\ge 1}  \frac{u^{n-1} \,z^n}{(1+tz)\cdots (n+tz)} \, \al_{n} \,.
\end{equation}
\end{Theorem}

As an immediate consequence of Theorem \ref{ThW}, we see
$c_{k,l}(u,t)=0$ unless $l> 0$.




\subsubsection{}

Equation \eqref{WAWc} is equivalent to
the equation
\begin{equation}
  \label{AA}
  \bA_k  =
\sum_{l\le k+1} c_{k,l}(u,t) \, \bA_0^l \,, \quad k=0,1,\dots\,,
\end{equation}
where the powers of $\bA_0$ are taken in the
associative algebra $\End(\infty)$.

The operator $\bA(z)$ is homogeneous of degree $-1$ with
respect to the following grading:
$$
\deg u = \deg t = - \deg z = 1 \,.
$$
Therefore, the operator $\bA_k$ has degree $k$ with respect
to the grading. Therefore, by \eqref{AA}, we see
\begin{equation}\label{homckl}
\deg\  c_{k,l}(u,t)\ = k \,.
\end{equation}
Theorem \ref{ThW} implies  $c_{k,l}(u,t)$
is a {\em monomial}:
\begin{equation}
  \label{cklt}
   c_{k,l}(u,t)= c_{k,l} \, u^{l-1} \, t^{k-l+1}\,, \quad  c_{k,l} \in \Q \,, 
\end{equation}
a nontrivial fact which will play an important role in the
proof.

Because of the homogeneity property
\eqref{homckl}, we may set $u=1$ in order to
 simplify
our computations.

\subsubsection{}

Taking the adjoint of equation \eqref{WAW} and reversing
the sign of $t$, we find: 
\begin{equation}
\label{WAWs}
\sum_{k\ge 0}  z^{k+1}\, \bAt_k^\star  =
\sum_{n\ge 1}  \frac{u^{n-1} \,z^n}{(1-tz)\cdots (n-tz)} \, \al_{-n} \,.
\end{equation}
where
$$
\bAt_k^\star = W^\star\, \bA^\star_k \, \left(W^\star\right)^{-1} \, ,
$$
and $W^\star= W^*(u,-t).$

Following the discussion of Section \ref{fullT}, 
we immediately obtain the following result.

\begin{Theorem}\label{TT2}
The triangular
linear change of time variables given by \eqref{WAW} and \eqref{WAWs}
makes the sequence of functions,
$$
\frac{q^{n^2/2}}{u^{n^2}} \,
e^{n u \partial} \, \tau(x,x^\star,u)\,, \quad  n\in \Z\,,
$$
a sequence of $\tau$-functions for
the full 2--Toda hierarchy of Ueno and Takasaki.
\end{Theorem}

Our derivation has neglected a minor point:
the operators $\bA_k, \bAt_k$
have constant terms when acting on $\LV$ (and similarly for 
$\bA_k^\star, \bAt_k^\star$).
However, these constants can be removed by
further conjugation by operators $\al_n$  in $\LV$.
The constants do not affect Theorem \ref{TT2}.

The explicit form of the linear change of variables from the
Gromov-Witten times to the standard times of the 2--Toda
hierarchy was conjectured by Getzler, see \cite{G2}.

\subsubsection{}
We now proceed with the proof of Theorem \ref{ThW} starting 
with the following result.

\begin{Proposition}\label{pThW1}
  For $k\ge 0$ and $l> 0$, the coefficient $c_{k,l}(u,t)$ is a monomial
in $t$ of degree $k-l+1$.
\end{Proposition}

\begin{proof}
We set $u=1$.
By \eqref{AA}, we may equivalently prove the
coefficient of $\bA_0^l$ in the expansion of
$\bA_k$ is a monomial in $t$ of degree
$k-l+1$. Further, by induction, it suffices to
prove the coefficients $b_{k,l}(t)$
in the expansion
\begin{equation}\label{AoAk}
\bA_0 \, \bA_{k} \, = \sum_{l \le k+1} b_{k,l}(t) \, \bA_l
\end{equation}
are monomials in $t$  of degree $k+1-l$ for $l\ge 0$.


The coefficients $b_{k,l}(t)$ with $l\ge 0$ can be determined
from the negative energy matrix elements of the product
$\bA_0 \, \bA_{k}$. The matrix elements of
$\bA_0 \, \bA_{k}$ are obtained as the $z\,w^{k+1}$
coefficient of the expansion of $\bA(z) \, \bA(w)$.
Since
$$
\cE_{a}(z) \, \cE_{b}(w) = e^{(aw-bz)/2} \, \cE_{a+b}(z+w) \,,
$$
we compute
\begin{multline}\label{AzAw}
\bA(z) \, \bA(w) \, = \cS(z)^{tz} \, \cS(w)^{tw} \,
\sum_{m\in\Z} \cE_{m}(z+w) \times \\
\frac{\cs(z)^{m} \, e^{mw/2}}{(1+tz)_m}\,
\left(\sum_{n\in\Z} \frac{(-tz-m)_{n}}{(1+tw)_n} \,
\left(\frac{1-e^{-w}}{1-e^z}\right)^n \right) \,.
\end{multline}
The summation over $n$ in \eqref{AzAw} is formally infinite,
but only finitely many terms actually contribute to
the $z\,w^{k+1}$ coefficient. Indeed, the coefficient of
$z$ vanishes if $m>n+1$, while the coefficient of $w^{k+1}$
vanishes if $n>k+1$.

The sum over $n$ in \eqref{AzAw} can be written as:
\begin{equation}
  \label{Fn}
  \Ff{-tz-m}{1+tw}{\frac{1-e^{-w}}{1-e^z}} +
\Ff{-tw}{1+m+tz}{\frac{1-e^z}{1-e^{-w}}} - 1 \, ,
\end{equation}
where the hypergeometric function is defined by \eqref{fS}.
The two series in \eqref{Fn} converge for
$|w|<|z|\ll 1$ and $|z|<|w|\ll 1$, respectively.
Therefore, we can write the coefficient of $z\,w^{k+1}$
as a sum of two contour integrals in two different
domains.

We may now deform these contour integrals to
integrals over
$$
|z|=|w|=\epsilon\ll 1 \,.
$$
The condition $m>0$
is needed for the hypergeometric function to remain
continuous in this limit. On the new contour, which is now common to
both integrals, we can use formula \eqref{lF1}.
After some simplifications, we find:
\begin{multline}\label{AAint}
  \bA_0 \, \bA_k = \frac{1}{(2\pi i)^2} \iint_{|z|=|w|=\varepsilon}
  \frac{dz\,dw}{z^2\,w^{k+2}} \, \times \\
  \frac{\left(1+\frac{w}{z}\right)^{tz+tw}}
  {\left(\frac{w}{z}\right)^{tw}}\,
  \frac{\Gamma(1+tz)\, \Gamma(1+tw)}{\Gamma(1+tz+tw)} \, \bA(z+w) + \dots
  \, 
\end{multline}
where the dots denote terms of non-negative energy.

The meaning of formula \eqref{AAint} is the following. First, the
multivalued function
\begin{equation} \label{mfun}
\frac{\left(1+\frac{w}{z}\right)^{tz+tw}}
  {\left(\frac{w}{z}\right)^{tw}}
\end{equation}
is defined using the cut
$$
\frac wz \ne (-\infty,0]\,.
$$
Because both $z$ and $w$ are small, the function \eqref{mfun} 
is integrable
in the neighborhood of the singularity $w = -z$ on the
contour of integration. Second, the negative energy terms
in $\bA(z+w)$ are nonsingular at $z+w=0$ and, hence, their
expansion in powers of $z$ and $w$ is
unambiguous. Also, these terms do not spoil the
convergence of the integral at $w = -z$.

{}From formula \eqref{AAint}, we deduce, for $l\ge 0$,
\begin{multline}\label{bklint}
 b_{k,l}(t) = \frac{1}{(2\pi i)^2} \sum_{a=0}^{l+1} \binom{l+1}{a}
\iint_{|z|=|w|=\varepsilon}
  \frac{dz\,dw}{z^{2-a}\,w^{k+a+1-l}} \, \times \\
  \frac{\left(1+\frac{w}{z}\right)^{tz+tw}}
  {\left(\frac{w}{z}\right)^{tw}}\,
  \frac{\Gamma(1+tz)\, \Gamma(1+tw)}{\Gamma(1+tz+tw)}
  \,.
\end{multline}
After replacing $tz$ and $tw$ by new variables, we see
\eqref{bklint} is indeed a monomial in $t$ of degree $k-l+1$.
\end{proof}

\subsubsection{}
{}From Lemma \ref{conjF}, we expect the following
heuristic result:
$$
\bA(z) \, \bA(w) \,\, \textup{``$=$''} \, \,
\frac{(z+w)^{tz+tw}}{z^{tz} \, w^{tw}} \,
\frac{(tz)! \, (tw)!}{(tz+tw)!} \, \bA(z+w) \,,
$$
which becomes a true equality when both $tz$ and $tw$ are
positive integers.
Equation \eqref{AAint} is a way to make sense of
of the heuristic formula.

\subsubsection{}
The next step in the proof of Theorem \ref{ThW} is the
following result.

\begin{Proposition}\label{PrW2}
For all $l$, the coefficient $c_{k,l}(u,t)$ is a monomial
in $t$ of degree $k-l+1$.
\end{Proposition}

\begin{proof}
By Proposition \ref{pThW1}, we need only consider
$l\le 0$. Define the operator $\bD$ by:
$$
\bD = W^{-1}  \left(\al_1 + H\right) W \,.
$$
Equation \eqref{divA} implies:
\begin{equation}
  \label{dvA}
  \left[\bD,\bAt_k\right] = t \, \bAt_{k-1} \,.
\end{equation}
Also, since $\bA_0=\al_1 + H + \dots$, we see
$$
\bD=\al_1 + \dots,
$$
where, in both cases, the dots stand for terms with
positive energy.

Since the matrix $\left[\bD,\al_1\right] = t\, \bAt_{-1}$
commutes with $\al_1$, the matrix $\bD$ has the
form
$$
\bD = \al_1 + \sum_{n>0} d_n(u,t) \, \al_{-n} \, H + \dots \,,
$$
where the dots stand for terms that
commute with $\al_1$ and
whose precise form depends on the ambiguity in the
choice of the dressing matrix $W$.
Here, $H$ is the energy operator and the
product is taken in
the algebra $\End(V)$.

It is easy to see equation \eqref{divA} uniquely
determines all the coefficients $d_n$ in terms of $c_{k,l}(u,t)$
with $l>0$. The coefficients $d_n$, in turn, determine all remaining
coefficients $c_{k,l}(u,t)$. In fact, 
$$
d_n(u,t) = - \frac{t^n}{u^n} \,.
$$
However, for the proof of the Proposition, we need only observe the
uniqueness forces
$d_n$ to have degree $n$ in $t$. Then,
the coefficients $c_{k,l}(u,t)$ must have degrees $k-l+1$ in $t$.
\end{proof}

\subsubsection{}

{}From the proof of Proposition \ref{PrW2}, we see
the matrices $\bAt_k$ can be uniquely characterized
by the two following conditions:
\begin{enumerate}
\item[(i)] $\bAt_0=\al_1$ and $\bAt_k$ is a linear combination
of $\al_1,\dots,\al_{k+1}$.
\item[(ii)] There exists a matrix of the form
$$
\bD=\al_{1}+\sum_{n>0} d_n \, \al_{-n} \, H\,, \quad d_1=-\frac{t}{u}\,,
$$
such that $\left[\bD,\bAt_k\right]=t\, \bA_{k-1}$ for
$k>0$.
\end{enumerate}

\subsubsection{}\label{endend}

We can now complete the proof of Theorem \ref{ThW}.
Since the coefficients $c_{k,l}(u,t)$ are monomials in $t$, the coefficients
are
identical to their leading order asymptotics as
$u\to 0$. Hence, the operators $\bAt_k$
can be determined by studying the $u\to 0$ asymptotics
of the operators $\bA(z)$. In the $u\to 0$ limit,
we have
\begin{multline}\label{bAasy}
  \bA(z) \sim
\sum_{n\ge 0} \frac{u^{n-1} \,z^n}{(1+tz)\cdots (n+tz)} \, \al_{n} + \\
\sum_{n>0} \frac{t}{u^{n+1}} \left(t-\frac 1z\right) \cdots
\left(t-\frac{n-1}z\right) \, \al_{-n} \,.
\end{multline}
In the $u\to 0$ limit, the dressing matrix $W$
becomes trivial and the statement of Theorem \ref{ThW}
can be read off directly from \eqref{bAasy}.

Formula \eqref{bAasy} also contains the
description of the dressed operators $\bAt_k$ for $k<0$.

\section{Commutation relations for operators $\A$}
\label{ffff}

Our goal here is to prove Theorem \ref{ThA}:
\begin{equation*}
   [\A(z,uz),\A(w,uw)]  = zw\, \de(z,-w) \,,
 \end{equation*}

\subsection{Formula for the commutators}
\label{aaa}

\subsubsection{}

We may calculate $\left[\A(z,uz),\A(w,uw)\right]$
by  the commutation relation (\ref{commcE}). We find,
\begin{equation}\label{commA}
[\A(z,uz),\A(w,uw)]  =
\cS(uz)^z \, \cS(uw)^w \,
\sum_{m\in\Z} c_{m}(z,w) \, \cE_{m}(u(z+w))
\end{equation}
where the functions $c_m(z,w)$ are defined by:
$$
c_m(z,w)=
\begin{cases}
\dfrac{\cs(u z)^{s} \, \cs(u w)^{s}}{(z+1)_{s}
(w+1)_{s}} \big[f_{s,u}(z,w)- f_{s,u}(w,z)\big],  & m=2s\,, \\
{} \\
\dfrac{\cs(u z)^{s} \, \cs(uw)^{s}}{(z+1)_{s}
(w+1)_{s}} \big[g_{s,u}(z,w)- g_{s,u}(w,z)\big],  & m=2s-1 \,.
\end{cases}
$$
Here, $f_{s,u}(z,w)$ and $g_{s,u}(z,w)$ are hypergeometric
series which are explicitly defined below.

We recall the definition of the hypergeometric series
which we require:
\begin{equation}
\Ff{-\nu}{\mu+1}{z} = \sum_{k=0}^\infty
\frac{\nu(\nu-1)\cdots(\nu-k+1)}{(\mu+1)\cdots(\mu+k)} \, (-z)^k \,,
\quad |z|<1\,.
\label{fS}
\end{equation}
Define $f_{s,u}(\mu,\nu)$ and $g_{s,u}(\mu,\nu)$ by:
\begin{multline}\label{lF2}
 f_{s,u}(\mu,\nu)= e^{-su\mu} \Ff{-\nu-s}{\mu+1+s}{\frac{1-e^{u\mu}}{1-e^{-u\nu}}} - \\
  e^{-su\nu} \Ff{-\nu-s}{\mu+1+s}{\frac{1-e^{-u\mu}}{1-e^{u\nu}}} +
  \frac{e^{-su\nu} - e^{-su\mu}}{2}
\end{multline}
and,
\begin{multline}\label{gu}
 g_{s,u}(\mu,\nu)= \frac{\nu+s}{\cs(u\nu)}\left[e^{(1-s)u\mu} \Ff{-\nu-s+1}{\mu+1+s}{\frac{1-e^{u\mu}}{1-e^{-u\nu}}} -\right. \\
\left.
  e^{-su\nu} \Ff{-\nu-s+1}
{\mu+1+s}{\frac{1-e^{-u\mu}}{1-e^{u\nu}}} \right] \,.
\end{multline}

\subsubsection{}

The series $f_{s,u}(z,w)$ and
$g_{s,u}(z,w)$ in formula \eqref{commA} are to be expanded in the ring
$\Q[u^{\pm}]((w))((z))$, that is, expanded in
Laurent series of $z$
with coefficients given by Laurent series in $w$.
Since, for example, the $k$th term in
\begin{multline}\label{sf}
\Ff{-w-s}{z+s+1}{\frac{1-e^{uz}}{1-e^{-uw}}} = \\
\sum_{k=0}^\infty
\frac{(w+s)\cdots(w+s-k+1)}{(z+s+1)\cdots(z+s+k)} \,
\left(
\frac{e^{uz}-1}{1-e^{-uw}}\right)^k
\end{multline}
is of order $z^k$, the extraction of any given
term in these expansions is, in principle, a finite
computation.
Similarly,
the series
 $f_{s,u}(w,z)$ and $g_{s,u}(w,z)$ in formula (\ref{commA}) are to be expanded
 in the ring
$\Q[u^{\pm}]((z))((w))$.

\subsubsection{}

The constant term of $\cE_0(u(z+w))$ plays a special role in formula
(\ref{commA}).
The expansion rules for the constant term,
\begin{equation}
\label{lpl}
\frac{f_{0,u}(z,w)}{\cs(u(z+w))} -
\frac{f_{0,u}(w,z)}{\cs(u(z+w))} \, ,
\end{equation}
are the following.
The first summand is to be expanded in ascending powers
of $z$ whereas the second summand is to be expanded in
ascending powers of $w$.

\subsubsection{}

We will show
the expansions of the two terms of $c_m(z,w)$ exactly
cancel each
other.
The commutator is therefore obtained entirely from the constant term.
We will show the expansions of the two terms of (\ref{lpl})
cancel except for
the two different expansions of the simple pole
at $z+w=0$.

\subsection{Some properties of the hypergeometric series}

\subsubsection{}

To proceed, several
properties of the hypergeometric
series \eqref{fS} are required.
Define the analytic continuation of \eqref{fS}
to the complex plane with
a cut along $[1,+\infty)$ by the following
integral:
\begin{equation}
  \label{fI}
  \Ff{-\nu}{\mu+1}z = \mu \int_0^1 (1-x)^{\mu-1} (1-zx)^{\nu} \, dx \,,\
\quad \Re \mu >0 \,.
\end{equation}
The above hypergeometric function is degenerate since
the elementary function,
$$
z^{-\mu} \, (1-z)^{\mu+\nu} \,,
$$
is a second solution to the hypergeometric equation and,
in addition, is  an eigenfunction of monodromy at
$\{0,1,\infty\}$.
As a consequence, the analytic continuation of the function
\eqref{fI} through the cuts $[1,+\infty)$
leads only to the appearance of elementary terms.
In fact, the analytic continuation of \eqref{fI}
through the cut $[1,+\infty)$ is given explicitly
by the formula \eqref{lF1} below.

\subsubsection{}

\begin{Lemma}\label{Lz1} For $z\notin [0,+\infty)$ we have:
\begin{multline}\label{lF1}
  \Ff{-\nu}{\mu+1}z = \\ 1- \Ff{-\mu}{\nu+1}{\frac1z} +
  \frac{(1-z)^{\mu+\nu}}{(-z)^\mu} \, \frac{\Gamma(\mu+1)\,
    \Gamma(\nu+1)}{\Gamma(\mu+\nu+1)} \,.
\end{multline}
\end{Lemma}

Here and in what follows we use the principal branches of
the functions $\ln w$ and $w^a$ for $w\notin (-\infty,0]$.

\begin{proof}
Integrating by parts and setting $y=zx$, the integral
\eqref{fI} is transformed to the following form:
$$
1-\nu \int_0^1 \left(1-\frac{y}{z}\right)^{\mu} (1-y)^{\nu-1} \, dy
+ \nu \int_z^1 \left(1-\frac{y}{z}\right)^{\mu} (1-y)^{\nu-1} \, dy\,.
$$
The last integral here is a standard beta-function integral and, thus,
the three terms in the above formula correspond precisely
to the three terms on the right side of \eqref{lF1} \,.
\end{proof}

\subsubsection{}

A similar argument proves the following result.

\begin{Lemma}\label{Lz2} For $z\notin [0,+\infty)$ we have
\begin{multline*}
  \nu\, \Ff{-\nu+1}{\mu+1}z = \\ \frac{\mu}{z}\,\,
 \Ff{-\mu+1}{\nu+1}{\frac1z} +
  \frac{(1-z)^{\mu+\nu-1}}{(-z)^\mu} \, \frac{\Gamma(\mu+1)\,
    \Gamma(\nu+1)}{\Gamma(\mu+\nu)} \,.
\end{multline*}
\end{Lemma}

\subsection{Conclusion of the proof of Theorem \ref{ThA}}

\subsubsection{}

\label{acc}
\begin{Lemma}\label{LF2}
The functions $f_{s,u}(\mu,\nu)$ and $g_{s,u}(\mu,\nu)$ are
analytic in a neighborhood of the origin
$(\mu,\nu)=(0,0)$ and symmetric in $\mu$ and $\nu$\,.
\end{Lemma}

\begin{proof}
We will prove the Lemma for $f_{s,u}(\mu,\nu)$. The
argument for $g_{s,u}(\mu,\nu)$ is parallel with Lemma
\ref{Lz2} replacing Lemma \ref{Lz1}.
The proof will show the
neighborhood can be
chosen to be independent of the parameter $s$.

For simplicity, we will first assume
$s$ is not a negative integer.
The assumption will be removed at the end of the proof.
Using relation \eqref{lF1}, we find,
\begin{equation}
f_{s,u}(\mu,\nu)=f_{s,u}(\nu,\mu)\label{ssff}
\end{equation}
on the intersection of the domains of
applicability of \eqref{lF1}.

The possible singularities of $f_{s,u}(\mu,\nu)$
near the origin are at $\nu=0$ and $\mu+\nu=0$,
corresponding to the singularities $z=\infty$ and
$z=1$ of the hypergeometric function (\ref{fI}), respectively.
The hypergeometric
function is analytic and single-valued in the complex
plane with a cut from $1$ to $\infty$.
The function $f_{s,u}(\mu,\nu)$ is well-defined if the arguments,
\begin{equation*}
\frac{1-e^{u\mu}}{1-e^{-u\nu}}\,, \frac{1-e^{-u\mu}}{1-e^{u\nu}}
\approx - \frac{\mu}{\nu} \, ,
\end{equation*}
do not fall on the cut $[1, + \infty)$.
Similarly, the function
$f_{s,u}(\nu,\mu)$ is well-defined
if the arguments,
$$\frac{1-e^{u\nu}}{1-e^{-u\mu}}\,, \frac{1-e^{-u\nu}}{1-e^{u\mu}}
\approx - \frac{\nu}{\mu} \, ,$$
do not fall on the cut $[1, +\infty)$.
By \eqref{ssff}, the
two functions above agree on the region where both are defined.
It follows
that $f_{s,u}(\mu,\nu)$
is single-valued and
analytic near the origin in the complement of
the divisor $\mu+\nu=0$. By
Lemma \ref{Lmv} below,  $f_{s,u}(\mu,\nu)$
remains bounded as $\nu\to -\mu$ and hence the
singularity at $\mu+\nu=0$ is removable.
We conclude $f_{s,u}(\mu,\nu)$ is analytic and symmetric
near the origin.

Finally, consider the case when $s\to -n$, where $n$ is
positive integer. The apparent simple pole
of $f_{s,u}(\mu,\nu)$ at $\mu= -s -n$ is, in fact, removable.
The removability follows either from symmetry (because there
is no such singularity in $\nu$) or else can be
checked directly using the formula
$$
\Res_{\mu=-n} \, \Ff{-\nu}{\mu}z
= (-1)^{n-1} \, \frac{(-\nu)_n}{(n-1)!} \, z^n \, (1-z)^{\nu-n} \,.
$$
\end{proof}

\subsubsection{}

\begin{Lemma}\label{Lmv}
We have
\begin{equation}\label{fs1}
f_{s,u}(\mu,-\mu)= - \frac{\mu}{s} \, \sinh(us\mu)
\end{equation}
and, in particular,
\begin{equation}
f_{0,u}(\mu,-\mu)= - u \mu^2\,.
\label{fs01}
\end{equation}
Similarly,
$$
g_{s,u}(\mu,-\mu)= \frac{s^2-\mu^2}{2s-1}\,
\frac{\sinh\frac{(2s-1)u\mu}2}
{\sinh\frac{u\mu}2} \,.
$$
\end{Lemma}

\begin{proof}
For $\Re s>0$ we can use the formula
$$
\Ff{\mu-s}{\mu+1+s}{1} =
\frac{\Gamma(\mu+1+s) \, \Gamma(2s)}{\Gamma(\mu+s) \Gamma(2s+1)} =
\frac{\mu+s}{2s} \,,
$$
from which \eqref{fs1} follows.
 By analytic continuation, \eqref{fs1}
holds for all $s$.  The computation of
 $g_{s,u}(\mu,-\mu)$ is identical.
\end{proof}

\subsubsection{}

We may now complete the proof of Theorem \ref{ThA}.
Since the functions $f_{s,u}(\mu,\nu)$ and
$g_{s,u}(\mu,\nu)$ are analytic
near the origin and symmetric in $\mu$ and
$\nu$, the nonconstant terms of formula \eqref{commA}
cancel.

The summands of the
constant term \eqref{lpl}  of formula \eqref{commA} can be analyzed using
\eqref{fs01}:
$$
\frac{f_{0,u}(z,w)}{\cs(u(z+w))} =
\frac{zw}{z+w} + \dots \,,
$$
where the dots represent a function analytic at
the origin and symmetric in $z$ and $w$.
Observe the prefactor in formula \eqref{commA}
is identically equal to 1 on the divisor
$z+w=0$ and does not affect the singularity.
The proof of Theorem \ref{ThA} is complete.

\vspace{+10 pt}
\noindent
Department of Mathematics \\
Princeton University \\
Princeton, NJ 08544\\
okounkov@math.princeton.edu \\

\vspace{+10 pt}
\noindent
Department of Mathematics\\
Princeton University\\
Princeton, NJ 08544\\
rahulp@math.princeton.edu
\end{document}